\documentclass[11pt]{article}
\usepackage{amsmath,amssymb,amsthm,mathrsfs}
\renewcommand{\mathcal}{\mathscr}

\newcommand{\cL}{\mathcal{L}}
\newcommand{\cT}{\mathcal{T}}
\newcommand{\cX}{\mathcal{X}}

\newcommand {\IR}{\mathbb{R}}

\newcommand {\F}{\mbox{${\mathcal F}$}}

\newcommand{\s}{\sigma}
\newcommand{\bm}{\mathbf}

\renewcommand{\P}{\mathrm{P}}
\newcommand {\E}{{\mathrm E}}
\newcommand{\1}{{\bf 1}}
\newcommand{\dimh}{\dim_{_{\rm H}}}
\renewcommand{\Cap}{{\rm Cap}}

\newcommand{\R}{\mathbb{R}}
\newcommand{\e}{\epsilon}

\newtheorem{stat}{Statement}[section]
\newtheorem{proposition}[stat]{Proposition}
\newtheorem{corollary}[stat]{Corollary}
\newtheorem{theorem}[stat]{Theorem}
\newtheorem{lemma}[stat]{Lemma}

\theoremstyle{definition} \newtheorem{remark}[stat]{Remark}
\newtheorem{defi}[stat]{Definition}

\newtheorem{nota}[stat]{Notation}

\numberwithin{equation}{section}

\begin{document}

\title{\bf Hitting probabilities for systems of 
	non-linear stochastic heat equations with additive noise}
\author{Robert C. Dalang$^{1,4}$, Davar Khoshnevisan$^{2,5}$, and Eulalia Nualart$^3$}
\date{}

\footnotetext[1]{Institut de Math\'ematiques, Ecole Polytechnique
	F\'ed\'erale de Lausanne, Station 8, CH-1015 
	Lausanne, Switzerland. \texttt{robert.dalang@epfl.ch}}

\footnotetext[2]{Department of Mathematics, The 
	University of Utah, 155 S. 1400 E. Salt Lake City, 
	UT 84112-0090, USA. \texttt{davar@math.utah.edu}}

\footnotetext[3]{Institut Galil\'ee, Universit\'e
	Paris 13, 93430 Villetaneuse, France.
	\texttt{nualart@math.univ-paris13.fr}}

\footnotetext[4]{Supported in part by the Swiss National Foundation for Scientific Research.}%

\footnotetext[5]{Research supported in part by a grant
	from the US National Science Foundation.}%
\maketitle
\begin{abstract}
	We consider a system of $d$ coupled non-linear stochastic
	heat equations in spatial dimension $1$ driven by $d$-dimensional additive 
	space-time white noise. We establish
	upper and lower bounds on hitting
	probabilities of the solution $\{u(t\,,x)\}_{t \in 
	\mathbb{R}_+, x \in [0\,,1]}$, in terms of respectively
	Hausdorff measure and Newtonian capacity.
	We also obtain the Hausdorff dimensions of level
	sets and their projections. A result of independent
	interest is an anisotropic form of the Kolmogorov continuity theorem.
\end{abstract}

 \vfill\pagebreak

\section{Introduction}

Let $\dot{W}:=(\dot{W}^1\,,\ldots,\dot{W}^d)$ be a vector of
$d$ independent space-time white noises
on $[0\,,T]\times[0\,,1]$. For all $1\le i\le d$, let
$b_i:\R^d\to\R$ be globally Lipschitz and
bounded functions, and ${\s}:=(\s_{i,j})$ be a
deterministic $d \times d$ invertible matrix (ellipticity).

Consider the system of stochastic partial differential equations
({\em s.p.d.e.}'s)
\begin{equation}\label{eq:spde}
    \frac{\partial u_i}{\partial t}(t\,,x) =
    \frac{\partial^2 u_i}{\partial x^2}(t\,,x)
    + \sum_{j=1}^d \s_{i,j} \dot{W}^j(t\,,x) +
    b_i(u(t\,,x)),
\end{equation}
for $1\le i\le d$, $t\in[0\,,T]$, and $x\in[0\,,1]$, where
${u}:=(u_1\,,\ldots,u_d)$, with initial conditions ${u}(0\,,x)={0}$
for all $x\in[0\,,1]$, and Neumann boundary conditions
\begin{equation}\label{neumann}
    \frac{\partial u_i}{\partial x}(t\,,0)=
    \frac{\partial u_i}{\partial x}(t\,,1) =0,
    \qquad 0\le t\le T.
\end{equation}

Equation (\ref{eq:spde}) is formal, but can be interpreted
rigorously  as follows (\textsc{Walsh} \cite{Walsh:86}):
Let $W^i=(W^i(s\,,x))_{s \in \mathbb{R}_+, \, x \in [0\,,1]}$,
$i=1,\ldots,d$, be independent Brownian sheets, defined on a
probability space $(\Omega\,,\F,\P)$, and
set $W=(W^1,\ldots,W^d)$. For $t \geq 0$, let
$\mathcal{F}_t=\sigma\{W(s\,,x),\ s \in [0\,,T],\ x \in [0\,,1]\}$.
We say that a process $u=\{u(t\,,x), \, t \in [0\,,T], \, x
\in [0\,,1]\}$ is {\em adapted} to $(\mathcal{F}_t)$ if $u(t\,,x)$ is
${\mathcal{F}}_t$-measurable for each $(t\,,x) \in [0\,,T] \times [0\,,1]$.
We say that $u = (u_1\,,\ldots,u_d)$ {\em is a solution of}
(\ref{eq:spde}) if $u$ is adapted to $(\mathcal{F}_t)$
and if for $i\in\{1\,,\ldots,d\}$,
$t\in[0\,,T]$, and $x\in[0\,,1]$,
\begin{equation}\begin{split} \label{spde}
    u_i(t\,,x) & = \int_0^t \int_0^1 G_{t-r}(x\,,v) \sum_{j=1}^d
        \s_{i,j}\, W^j(drdv)
        + \int_0^t\int_0^1 G_{t-r}
        (x\,,v)\, b_i(u(r\,,v)) \, dr\, dv.
\end{split}\end{equation}
Here, $G_t(x\,,y)$ denotes the Green kernel for the heat equation
with Neumann boundary conditions. See, for example,
\textsc{Walsh} \cite{Walsh:86} or \textsc{Bally, Millet, and Sanz--Sol\'e} \cite{Bally:95}.

Our goal is to develop aspects of
potential theory for the solution
to the system of stochastic heat equations \eqref{eq:spde}.
In particular, given $A \subset \IR^d$, we want to determine
whether the process 
$\{u(t\,,x), \, t \geq 0, \, x \in [0\,,1]\}$ visits, or hits, $A$
with positive probability.

Potential theory for single-parameter processes is a mature subject.
See, for example \textsc{Blumenthal and Getoor}
\cite{Blumenthal:68}, \textsc{Port and Stone}
\cite{Port:78}, and \textsc{Doob} \cite{Doob:84}.
There is also a growing literature on the
potential theory for multiparameter processes
(\textsc{Khoshnevisan} \cite{Khoshnevisan:02}.

For the linear form of (\ref{eq:spde})
$(b \equiv 0, \sigma \equiv I_d$, where $I_d$
denotes the $d \times d$ identity matrix), results on hitting probabilites 
have been obtained in \textsc{Mueller and Tribe} \cite{Mueller:03}. In the case $d=1$, 
for a particular form of (\ref{eq:spde}) with additive noise 
$(\sigma \equiv I_d, \, b(u) = u^{-\delta}$ for $\delta > 3$
and $b(u) = cu^{-3})$, the issue of whether or not the solution hits 
$0$ has been discussed in \textsc{Zambotti}
\cite{Zambotti:02,Zambotti:03} and 
\textsc{Dalang, Mueller, and Zambotti} \cite{Dalang:06}.

For {\em non-linear} s.p.d.e.'s, a general
result was obtained in \textsc{Dalang and Nualart}
\cite{Dalang:04}, valid for systems of reduced hyperbolic
equations on $\R_+^2$ (essentially equivalent to systems
of wave equations in spatial dimension 1) that are driven by
two-parameter white noise. In this paper, we will be concerned
with obtaining upper and lower bounds on hitting probabilities
for the solution of the system (\ref{eq:spde}). In a
forthcoming paper \cite{Dalang-K-N:07}, we use quite different
techniques from the Malliavin calculus, consider systems
of non-linear heat equations with multiplicative noise, and
obtain bounds that are slightly different than those in this paper.

Let $\{v(r)\}_{r\in T}$ denote a random field that takes
values in $\R^d$, where $T$ is some Borel-measurable
subset of $\R^N$. Let $v(T)$ denote the range of $T$
under the random map $r \mapsto v(r)$.
We say that a Borel set $A\subseteq\R^d$
is called \emph{polar} for $v$ if $\P\{v(T)\cap A\neq\varnothing\}=0$;
otherwise, $A$ is called \emph{nonpolar}.
Two of our main results are the following. They will be proved in
Section \ref{sec5}.

\begin{theorem}\label{th:hitprob}
	Let ${u}$ denote the solution to
	\textnormal{(\ref{eq:spde})} on $]0\,,T]\times [0\,,1]$.
	\begin{itemize}
		\item[\textnormal{(a)}] A (nonrandom) Borel set $A\subset\R^d$
			is nonpolar for $(t\,,x) \mapsto u(t\,,x)$ if
			it has positive $(d-6)$-dimensional
			capacity. On the other hand, if $A$ has zero
			$(d-6)$-dimensional Hausdorff measure, then
			$A$ is polar for $(t,x) \mapsto u(t\,,x)$.
		\item[\textnormal{(b)}] Fix $t\in\,]0\,,T]$. A Borel set
			$A\subseteq\R^d$ is nonpolar for $x\mapsto u(t\,,x)$ if
			$A$ has positive $(d-2)$-dimensional capacity. If, on the other
			hand, $A$ has zero $(d-2)$-dimensional Hausdorff measure, then
			$A$ is polar for $x\mapsto u(t\,,x)$.
		\item[\textnormal{(c)}] Fix $x\in [0\,,1]$. A Borel set
			$A\subseteq\R^d$ is nonpolar for $t\mapsto u(t\,,x)$ if
			$A$ has positive $(d-4)$-dimensional capacity. If, on the other
			hand, $A$ has zero $(d-4)$-dimensional Hausdorff measure, then
			$A$ is polar for $t\mapsto u(t\,,x)$.
	\end{itemize}
\end{theorem}

The definitions of capacity and Hausdorff measures will be recalled 
shortly.

There is a small gap between the conditions
of positive capacity and positive Hausdorff measure.
In some cases, we know
how to bridge that gap. Indeed, the results of
\textsc{Mueller and Tribe} \cite{Mueller:03} will make
this possible in parts (a) and (b) of the following.
This reference does not however apply to statement (c).

\begin{corollary}\label{cor:hitprob}
	Let ${ u}$ denote the solution to \textnormal{(\ref{eq:spde})}.
	\begin{itemize}
		\item[\textnormal{(a)}]
			Singletons are polar for 
			$(t\,,x) \mapsto { u}(t\,,x)$ if and only if $d \geq 6$.
		\item[\textnormal{(b)}] Fix $t\in\,]0\,,T]$. 
			Singletons are polar for $x \mapsto { u}(t\,,x)$
			if and only if $d \geq 2$.
		\item[\textnormal{(c)}]
			Fix $x\in [0\,,1]$. Singletons are polar for $t \mapsto { u}(t\,,x)$
			if $d >4$ and are nonpolar when $d < 4$.
			The case $d=4$ is open.
	\end{itemize}
\end{corollary}

\noindent This corollary is proved in Section \ref{sec5}.

Our work has other, ``more geometric,'' consequences as well.
For instance, in Corollary \ref{codim2} below, we prove that
if $d\ge 6$, then the Hausdorff dimension of the range of the solution
to \eqref{eq:spde} is $6$ a.s. On the other hand, when $d<6$,
Corollary \ref{cor:hitprob} implies readily that the range
of the solution to \eqref{eq:spde} has full Lebesgue measure a.s.

This paper is organized as follows. In Section \ref{sec2}
we present general conditions on an $\R^d$-valued random field
$(v(t\,,x))$ that imply {\em lower} bounds on hitting
probabilities (Theorem \ref{thm:cap:LB}). These conditions
are stated in terms of a lower bound
on the one-point density function of the
random vectors $v(t\,,x)$ and an upper bound on the two-point
density function; that is, the density function of
$(v(t\,,x), v(s\,,y))$
for $(t\,,x) \neq (s\,,y)$ (see conditions {\bf A1} and {\bf A2}).
These conditions also yield information about level sets of
the process and their projections (Theorem \ref{thm:cap2:LB}).
They are related to, but not identical with, 
the conditions of \textsc{Dalang and Nualart}
\cite{Dalang:04}.

In Section \ref{sec3} we isolate properties of the random
field that imply {\em upper} bounds on hitting probabilities
(Theorem \ref{thm:meas:UB}), and corresponding properties
of level sets and their projections (Theorem \ref{thm:cap2:LBB}).
These conditions are implied by sufficient conditions that are
often not too difficult to check, namely that the one-point density
function of the random variables $v(t,x)$ is uniformly bounded above
and an estimate on $L^p$-moments of increments of the random field
(Theorem \ref{eta}), similar to the condition in the classical
Kolmogorov continuity theorem. These conditions are different
from those of \cite{Dalang:04} which made specific use
of the structure of the filtration of the solution to a hyperbolic
s.p.d.e.\ in $\R_+^2$, and, in particular, of Cairoli's maximal
inequality for $2$-parameter martingales; there is no counterpart
to these for the stochastic heat equation.

In Section \ref{sec4} we verify the conditions of Sections
\ref{sec2} and \ref{sec3} for the solution of the linear form
of (\ref{eq:spde}), that is, with $b \equiv 0$ (see
Theorem \ref{rd:thm4.9}). In order to obtain the best
estimates possible, a careful analysis of moments of
increments and of the determinant of the variance/covariance
matrix of the (in this case, Gaussian) process $(u(t\,,x))$ is
needed. This also requires a version of the Kolmogorov
continuity theorem that is tailored to the needs of the
stochastic heat equation. This is presented in Appendix
\ref{apendB}, and may be of independent interest.

Finally, in Section \ref{sec5}, we use Girsanov's
theorem to transfer results about hitting probabilities
of the solution to the linear form of (\ref{eq:spde})
to the general form of (\ref{eq:spde}) (Proposition
\ref{final}), and we prove Theorem \ref{th:hitprob} and
Corollary \ref{cor:hitprob}. Some results on capacity
and energy are gathered in Appendix \ref{apendA}.

Let us conclude this Introduction by defining the requisite
notation and terminology. For all Borel sets $F\subseteq\R^d$
we define $\mathcal{P}(F)$ to be the set of all probability measures
with compact support in $F$.
For all integers $k\ge 1$ and
$\mu\in\mathcal{P}(\R^k)$, we let $I_\beta(\mu)$ denote the
\emph{$\beta$-dimensional energy} of $\mu$; that is,
\begin{equation}
    I_\beta(\mu)  := \iint {\rm K}_\beta(\|x-y\|)\, \mu(dx)\,\mu(dy),
\end{equation}
where $\|x\|$ denotes the Euclidean norm of $x \in \R^k$. Here
and throughout,
\begin{equation} \label{k}
	{\rm K}_\beta(r) := 
	\begin{cases}
		r^{-\beta}&\text{if $\beta >0$},\\
		\log  ( N_0/r ) &\text{if $\beta =0$},\\
		1&\text{if $\beta<0$},
    \end{cases}
\end{equation}
where $N_0$ is a constant whose value will be specified later
in the proof of Lemma \ref{prel}.

If $f: \R^d \mapsto \R_+$ is a probability density function,
then we will write $I_\beta(f)$
for the $\beta$-dimensional energy of the measure $f(x)dx$.

For all $\beta\in\R$, integers $k\ge 1$, and Borel sets
$F\subset\R^k$, $\text{Cap}_\beta(F)$ denotes the
\emph{$\beta$-dimensional capacity of
$F$}; that is,
\begin{equation}
    \text{Cap}_\beta(F) := \left[ \inf_{\mu\in\mathcal{P}(F)}
    I_\beta(\mu) \right]^{-1},
\end{equation}
where $1/\infty:=0$.

Given $\beta\geq 0$, the $\beta $-dimensional \emph{Hausdorff measure}
of $F$ is defined by
\begin{equation}\label{eq:HausdorfMeasure}
	{\mathcal{H}}_ \beta (F)= \lim_{\epsilon \rightarrow 0^+} \inf
	\left\{ \sum_{i=1}^{\infty} (2r_i)^ \beta : F \subseteq
	\bigcup_{i=1}^{\infty} B(x_i\,, r_i), \ \sup_{i\ge 1} r_i \leq
	\epsilon \right\},
\end{equation}
where $B(x\,,r)$ denotes the open (Euclidean) ball of radius
$r>0$ centered at $x\in \R^d$. When $\beta <0$, we define 
$\mathcal{H}_\beta (F)$ to be infinite.

Throughout, we consider the following \emph{parabolic metric}:
For all $s,t\in[0\,,T]$ and $x,y\in[0\,,1]$,
\begin{equation}\label{eq:Delta}
    {\bf \Delta}((t\,,x)\,; (s\,,y)) :=
    |t-s|^{1/2} + |x-y|.
\end{equation}
Clearly, this is a metric on $\R^2$ which
generates the usual Euclidean topology on $\R^2$.
We associate to this metric the energy form
\begin{equation}
    I_\beta^{{\bf \Delta}} (\mu) := \iint {\rm K}_\beta(
    {\bf \Delta}((t,x)\,;(s,y)))
    \,\mu(dt\,dx)\,\mu(ds\,dy),
\end{equation}
and its corresponding capacity
\begin{equation}
    \text{Cap}^{\bf \Delta}_\beta (F) := \left[
    \inf_{\mu\in\mathcal{P}(F)} I_\beta^{ \bf \Delta}(\mu)
    \right]^{-1}.
\end{equation}
For the Hausdorff measure, we write
\begin{equation}
	{\mathcal{H}}^{\bf \Delta}_s (F)= \lim_{\epsilon \rightarrow 0^+} \inf
	\biggl\{ \sum_{i=1}^{\infty} (2r_i)^s : F \subseteq
	\bigcup_{i=1}^{\infty} B_{\bf \Delta}((t_i\,,x_i)\,, r_i), \ 
	\sup_{i\ge 1} r_i \leq
	\epsilon \biggr\},
\end{equation}
where $B_{\bf\Delta}((t\,,x)\,,r)$ denotes the open
$\bf\Delta$-ball of radius $r>0$ centered at $(t\,,x)\in[0\,,T]\times[0\,,1].$

\section{Lower Bounds on Hitting Probabilities}\label{sec2}

Fix two compact intervals $I$ and $J$ of $\R$.
Suppose that $\{v(t\,,x)\}_{(t,x)\in I \times J}$ is a two-parameter, \emph{continuous}
random field with values in $\R^d$, such that $(v(t\,,x)\,,v(s\,,y))$  
has a joint probability density function $p_{t,x;s,y}(\cdot\,,\cdot)$,
for all $s,t\in I$ and
$x,y\in J$ such that $(t\,,x) \neq (s\,,y)$. That is,
\begin{equation}
    \E \left[ f\left( v(t\,,x)\,,  v(s,y)\right) \right]
    = \iint f(a\,,b) \, p_{t,x;s,y}(a\,,b)\, da\, db,
\end{equation}
for all bounded Borel-measurable functions $f:I \times J \to\R$.
We will denote the marginal density function of $v(t\,,x)$ by
$p_{t,x}$.

\vskip 12pt

Consider the following hypotheses:
\begin{description}
    \item[A1.] For all $M>0$, there exists a positive
    	and finite constant $C
        = C(I,J,M,d)$ such that for all
        $(t\,,x) \in I\times J$ and all $z\in[-M\,,M]^d$,
        \begin{equation}
             p_{t,x}(z) \ge C.
        \end{equation}
    \item[A2.] There exists $\beta>0$ such that for all
        $M>0$, there exists $c = c(I,J,\beta,M,d)>0$
        such that for all $s,t\in I$ and $x,y\in J$
        with $(t\,,x) \neq (s\,,y)$,
        and for every $z_1,z_2\in[-M\,,M]^d$,
        \begin{equation}
            p_{t,x;s,y}(z_1\,,z_2) \le
            \frac{c}{\big[ {\bf \Delta}((t\,,x)\,; (s\,,y))
            \big]^{\beta/2}}\exp\left(
            - \frac{\| z_1-z_2\|^2}{c
            {\bf \Delta}((t\,,x)\,; (s\,,y))}\right).
        \end{equation}
\end{description}

Our next theorem discusses lower bounds for various hitting probabilities
of the random field $v$.

\begin{theorem}\label{thm:cap:LB}
    Suppose {\bf A1} and {\bf A2} are met. Fix $M >0$.
    \begin{enumerate}
        \item[\textnormal{(1)}]  There exists a positive and finite constant $a=a(I,J,\beta,M,d)$ such that for all compact sets $A\subseteq[-M\,,M]^d$,
            \begin{equation}
                \P\left\{ v(I\times J) \cap A
                \neq\varnothing \right\} \ge a\Cap_{\beta-6}(A).
            \end{equation}
        \item[\textnormal{(2)}] There exists a positive and finite constant
            $a=a(J,M\,\beta,d)$ such that for all $t\in I$ and for all compact
            sets $A\subseteq[-M\,,M]^d$,
            \begin{equation}
                \P\left\{ v(\{t\}\times J) \cap A
                \neq\varnothing \right\} \ge a\Cap_{\beta-2}(A).
            \end{equation}
        \item[\textnormal{(3)}] There exists a positive and finite constant
            $a=a(I,M,\beta,d)$ such that for all $x \in J$ and for all compact
            sets $A\subseteq[-M\,,M]^d$,
            \begin{equation}
                \P\left\{ v(I\times \{x\}) \cap A
                \neq\varnothing \right\} \ge a\Cap_{\beta-4}(A).
            \end{equation}
    \end{enumerate}
\end{theorem}

Before proving this theorem, we need two technical lemmas.

\begin{lemma} \label{prel}
		Fix $N>0$ and $\beta>0$.
	    \begin{enumerate}
	    \item[\textnormal{(1)}] There exists a finite and positive constant
		    $C_1=C_1(I,J,\beta,N )$ such that for all $a\in[-N\,,N]$,
		    \begin{equation}
		        \int_{I} dt \int_{I} ds  \int_{J} d x \int_{J} dy
		        \
		        \frac{e^{-a^2/{\bf \Delta}((t,x);(s,y))}}{{\bf \Delta}^{\beta/2}
		        ((t,x);(s,y))}
		        \le C_1 \, {\rm K}_{\beta-6}(a).
		    \end{equation}
	    \item[\textnormal{(2)}] Fix $\alpha>0$. There exists a finite and positive constant
		    $C_2=C_2(I,\beta, N)$ such that for all $a\in[-N\,,N]$,
		    \begin{equation}
		        \int_{I} d t \int_{I} ds \
		        \frac{e^{-a^2/\vert t-s\vert^\alpha}}{\vert t-s\vert^{\alpha\beta/2}}
		        \le  C_2 \, {\rm K}_{\beta-(2/\alpha)}(a).
		    \end{equation}
	    \end{enumerate}
\end{lemma}

\begin{proof}
	   We start by proving (1). Using the change of variables
		$\tilde{u}=t-s$ ($t$ fixed), $\tilde{v}=x-y$ ($x$ fixed), we have
		\begin{equation}\begin{split}
			& \int_{I} dt \int_{I} ds  \int_{J} d x \int_{J} dy
		        \,
		        \frac{e^{-a^2/{\bf \Delta}((t,x);(s,y))}}{{\bf \Delta}^{\beta/2}
		        ((t,x);(s,y))}\\
			&\hskip1in \leq 4 \vert I \vert \, \vert J \vert 
				\int_0^{\vert I \vert} d\tilde{u}
				\int_0^{\vert J \vert} d\tilde{v} \, 
				(\tilde{u}^{1/2}+\tilde{v})^{-\beta/2} \exp \biggl(-
				\frac{a^2}{\tilde{u}^{1/2}+\tilde{v}} \biggr).
		\end{split}\end{equation}
		A change of variables [$\tilde{u}=a^4 u^2$, $\tilde{v}=a^2 v$] 
		implies that this is equal to
		\begin{equation} \label{chv}
			C a^{6-\beta} \int_0^{r} du  \int_0^{m} dv \, \frac{u}{(u+v
			)^{\beta/2}} \exp
			\biggl(-\frac{1}{u+v} \biggr),
		\end{equation}
		where $r:=\sqrt{\vert I \vert}/a^2$ and
		$m:=\vert J \vert/a^2$. Notice that 
		$r \geq r_1:= \sqrt{\vert I \vert}/N^2>0$ and
		$m \geq m_1:= \vert J \vert/N^2>0$.

		Observe that
		\begin{equation}
			\int_0^{r} du \int_0^{m} dv \, \frac{u}{(u+v )^{\beta/2}} \exp
			\biggl(-\frac{1}{u+v} \biggr) \leq \int_0^{r} du 
			\int_0^{m} dv \, (u+v )^{1-\frac{\beta}{2}} \exp
			\biggl(-\frac{1}{u+v} \biggr).
		\end{equation}
		Pass to polar coordinates to deduce that the preceding is 
		bounded above by $I_1+I_2(r\,,m)$, where
		\begin{equation}\begin{split}
			&I_1:= \int_0^{\sqrt{r^2_1+m_1^2}} d\rho \,
				\rho^{2-(\beta/2)} \exp (-c/\rho), \\
			&I_2(r\,,m):=\int_{\sqrt{r^2_1+m_1^2}}^{\sqrt{r^2+m^2}} d\rho \,
				\rho^{2-(\beta/2)} .
		\end{split}\end{equation}
		Clearly, $I_1 \leq C < \infty$, and if $\beta \neq 6$, then
		\begin{equation}
			I_2(r\,,m)=\frac{\left( \sqrt{r^2+m^2} \right)^{3-(\beta/2)}-
			\left( \sqrt{r^2_1+m^2_1} \right)^{%
			3-(\beta/2)}}{3-(\beta/2)}.
		\end{equation}
		There are three separate cases to consider:
		(i) If $\beta > 6$, then $3-(\beta/2)<0$, 
		and hence $I_2(r\,,m) \leq C$ for all 
		$r \geq r_1$ and $m \geq m_1$.
		(ii) If $\beta<6$, then $I_2(r\,,m) \leq c 
		(\sqrt{r^2+m^2})^{3-(\beta/2)}=C a^{\beta-6}$
		for all $r \geq r_1$ and $m \geq m_1$. 
		(iii) Finally, if $\beta=6$, then
		\begin{equation}\begin{split}
			I_2(r\,,m) &\leq C \left[ \ln \left(\sqrt{r^2+m^2} \right)-
				\ln(r_1) \right]\\
			&=c\left[ \ln\left(
				\frac{\vert I \vert+\vert J \vert^2}{r_1}\right)+
				2 \ln \left(\frac{1}{a} \right) \right].
		\end{split}\end{equation}
		
		We combine these observations to deduce that for all $\beta>0$
		there exists $C>0$ such that for all $a \in [-N\,,N]$,
		the expression in (\ref{chv}) is bounded above by
		\begin{equation}
			C a^{6-\beta} (I_1+I_2(r\,,m)) \leq c \, {\rm K}_{\beta-6}(a),
		\end{equation}
		provided that $N_0$ in (\ref{k}) is sufficiently large.
		This proves (1).

		Next we prove (2).  Fix $t$ and change variables
		[$u=t-s$] to see that
		\begin{equation}
			\int_{I} d t \int_{I} ds
			\
			\frac{e^{-a^2/\vert t-s\vert^\alpha}}
			{\vert t-s\vert^{\alpha \beta/2}}
			\leq 2 \vert I \vert
			\int_0^{\vert I \vert} du \, u^{-\alpha \beta/2}
			e^{-a^2/u^\alpha}.
		\end{equation}
		Another change of variables [$u=a^{2/\alpha} v$] simplifies this expression to
		\begin{equation} \label{chv2}
			C a^{(2/\alpha)-\beta} \int_0^{r} dv \, v^{-\alpha \beta/2} 
			e^{-1/v^{\alpha}},
		\end{equation}
		where $r:=\vert I \vert\, a^{-2/\alpha}$. Notice that 
		$r \geq r_1:= \vert I \vert\, N^{-2/\alpha}>0$.
		
		Observe that
		\begin{equation} 
			\int_0^{r} dv \, v^{-\alpha \beta/2} e^{-1/v^{\alpha}} \leq I_1+I_2(r),
			\end{equation}
			where
\begin{equation}
			I_1:= \int_0^{r_1} dv \, v^{-\alpha \beta/2} e^{-1/v^{\alpha}}, \qquad
			I_2(r):=\int_{r_1}^{r} dv \, v^{-\alpha \beta/2} .
\end{equation}
		Clearly, $I_1 \leq C <\infty $. Moreover, if $\alpha \beta \neq 2$ then 
		\begin{equation}
			I_2(r)=\frac{r^{1-(\alpha \beta/2)}-r_1^{1-(\alpha \beta/2)}}{1-(\alpha \beta/2)}.
		\end{equation}
As above, we consider three different cases: (i) If $\alpha \beta>2$, then $1-(\alpha \beta/2)<0$, and hence
		$I_2(r) \leq C$ for all $r \geq r_1$. (ii) If $\alpha \beta<2$, then $I_2(r) \leq C a^{-(2/\alpha)+\beta}$ for all $r \geq r_1$. 
		(iii) If $\alpha \beta =2$, then
			\begin{equation}
			I_2(r) =\left[ \ln\left(
			\frac{\vert I \vert}{r_1}\right)+
			\frac{2}{\alpha} \ln \left(\frac{1}{a} \right) \right].
		\end{equation}
		We combine these observations to deduce that for all $\beta>0$ and $\alpha>0$,
		there exists $C>0$ such that for all $a \in [-N\,,N]$,
		the expression in (\ref{chv2}) is bounded above by
		\begin{equation}
			C a^{(2/\alpha)-\beta} (I_1+I_2(r)) \leq c \, {\rm K}_{\beta-(2/\alpha)}(a),
		\end{equation}
		provided that $N_0$ in (\ref{k}) is sufficiently large.
		This proves (2) and completes the proof of the lemma.	
\end{proof}

For all $a,\nu,\rho>0$, define 
\begin{equation} \label{rdpsi}
    \Psi_{a,\nu}(\rho) := \int_0^a \frac{dx}{\rho + x^\nu}.
\end{equation}

\begin{lemma}\label{lem:Psi}
    For all $a,\nu, T>0$, there exists a finite
    and positive constant $C=C(a\,,\nu \,,T)$ such that
    for all $0<\rho<T$,
    \begin{equation}
        \Psi_{a,\nu}(\rho) \le C {\rm K}_{(\nu-1)/\nu}(\rho).
    \end{equation}
\end{lemma}

\begin{proof}
    If $\nu<1$, then $\lim_{\rho \to 0}
    \Psi_{a,\nu}(\rho)=\int_0^a x^{-\nu}\, dx<\infty$.
    In addition, $\rho \mapsto \Psi_{a,\nu}(\rho)$ is nonincreasing, so
$\Psi_{a,\nu}$ is bounded on $\R_+$
    when $\nu<1$. In this case, ${\rm K}_{(\nu-1)/\nu}(\rho)=1$,
    whence follows the result in the case that $\nu<1$.

    For the case $\nu\ge 1$, we change variables $(y=x\rho^{-1/\nu})$
    to find that
    \begin{equation}
        \Psi_{a,\nu}(\rho) = \rho^{-(\nu-1)/\nu}
        \int_0^{a \rho^{-1/\nu}} \frac{dy}{
        1+ y^\nu}.
    \end{equation}
    When $\nu>1$, this gives the desired result, with $c= \int_0^{+\infty} dy\, (1+ y^\nu)^{-1}$. 
    When $\nu=1$, we simply evaluate the integral in (\ref{rdpsi}) explicitly: this gives the result for $0< \rho < T$,
    given the choice of ${\rm K}_0(r)$ in (\ref{k}). 
    We note that the constraint ``$0 < \rho < T$'' is needed only in this case.
\end{proof}

\begin{proof}[Proof of Theorem \ref{thm:cap:LB}] 
	We begin by proving (1). Let $A \subset [-M,M]^d$ be a
	compact set. Without loss of generality, we assume
	that $\Cap_{\beta-6}(A)>0$; otherwise there is nothing to prove.
    By Taylor's theorem (cf. \textsc{Khoshnevisan}
    \cite[Appendix C, Corollary 2.3.1, p.525]{Khoshnevisan:02})
    this implies that $\beta -6 <d$ and $A \neq \varnothing$.

    There are separate cases to consider:
    \vskip 12pt

    {\em Case 1:}  $\beta-6 <0$. Then $\Cap_{\beta-6}(A)=1$.
    Hence it suffices to prove that there exists a finite 
    and positive constant $a$ (that does not depend on $A$) such that
    \begin{equation} \label{positive}
			\P\left\{ v(I\times J) \cap A
			\neq\varnothing \right\} \geq a.
		\end{equation}
		
      Define, for all $z \in \mathbb{R}^d$ and $\epsilon>0$, 
    $\tilde{B}(z \,,\e):=\{y\in\R^d:\ |y-z|<\e\}$, where $|z|:=\max_{1\le j\le d}|z_j|$, and
    \begin{equation}
        J_\e(z) = \frac{1}{(2\e)^d} \int_{I} dt \, \int_{J } dx \, \1_{\tilde{B}(z,\e)}
        (v(t\,,x)).
    \end{equation}
    Fix $z \in A\subseteq[-M\,,M]^d$. Hypothesis {\bf A1}
    implies that for all $\epsilon>0$,
	\begin{equation} \label{Jmin}\begin{split}
		\E \left[ J_\e(z) \right] &= \frac{1}{(2\e)^d}
			\int_{I} dt \, \int_{J } dx \, \int_{\tilde{B}(z,\e)}
			da\, p_{t,x}(a)\\
		&\ge C \vert I \vert\, \vert J \vert,
	\end{split}\end{equation}
    where $C>0$ does not depend on $z$.

    On the other hand, {\bf A2} implies that
\begin{equation}\begin{split}
        \E \left[ (J_\e(z))^2 \right] &= \frac{1}{(2\e)^{2d}}
        \int_I dt \, \int_J dx \, \int_I ds \, \int_J dy \,
        \int_{\tilde{B}(z,\e)} dz_1  \int_{\tilde{B}(z,\e)} dz_2
        \ p_{t,x;s,y}(z_1,z_2) \\ 
        & \leq c \int_I dt  \, \int_J dx\, \int_I ds \,
        \int_J dy \, \frac{1}{[ {\bf \Delta} ((t\,,x)\,; (s\,,y)) ]^{\beta/2}}.
\end{split}\end{equation}
    The change of variables
		$u=t-s$ ($t$ fixed), $v=x-y$ ($x$ fixed), implies that the preceding is bounded above by
		\begin{equation}
        C \int_0^{ \vert I \vert} du \,  \int_0^{\vert J \vert}
        dv \,(u^{1/2}+v)^{-\beta/2} \leq C' \int_0^{ \vert I \vert}
        du \, \Psi_{\vert J \vert, \beta/2} (u^{\beta/4}).
    \end{equation}
    Therefore, Lemma \ref{lem:Psi} implies that for all $\epsilon>0$,
    \begin{equation}
        \E \left[ (J_\e(z))^2 \right] \leq C \int_0^{ \vert I \vert}
        du \, {\rm K}_{1-(2/\beta)}(u^{\beta/4}).
    \end{equation}
    In order to bound the preceding integral, consider three
    different cases: (i) If $0<\beta<2$, then $1-{2/\beta}<0$
    and the integral equals $\vert I \vert$.
    (ii) If $2<\beta<6$, then ${\rm K}_{1-(2/\beta)}(u^{\beta/4})=
    u^{(1/2)-(\beta/4)}$ and the integral is finite. (iii) If $\beta=2$, then 
    ${\rm K}_{0}(u^{\beta/4})=\text{log}(N_0/u^{1/2})$ and the integral is also finite. 
    This fact, (\ref{Jmin}), and the Paley--Zygmund inequality
    (\textsc{Khoshnevisan} \cite[Lemma 1.4.1, Chap.3]{Khoshnevisan:02})
    together imply that
     \begin{equation}
    \P\left\{ J_\e(z) >0 \right\} \geq C >0.
		\end{equation} 
		 The left-hand side is bounded above by
    $\P\{ v(I \times J) \cap A^{(\e)}\neq\varnothing\}$,
    where $A^{(\e)}$ denotes the closed $\e$-enlargement
    of $A$. Let $\e\downarrow 0$ and appeal to the continuity
    of the trajectories of $v$ to find that 
    \begin{equation}
        \P\left\{ v(I\times J) \cap A \neq \varnothing\right\}
        \ge C>0.
    \end{equation}
   This proves (\ref{positive}).

   \vskip 12pt

    {\em Case 2:}  $0< \beta-6 < d$. Define, for all $\mu \in  \mathcal{P}(A)$ and
    $\e>0$,
    \begin{equation}
        J_\e(\mu) = \frac{1}{(2\e)^d}
        \int_{\R^d} \mu(dz) \, \int_{I} dt \, \int_{J} dx \, \1_{\tilde{B}(z,\e)}
        (v(t\,,x)).
    \end{equation}
    Fix $\mu \in \mathcal{P}(A)$ such that
    \begin{equation}\label{eq:mu1}
        I_{\beta-6}(\mu) \le \frac{2}{\Cap_{\beta-6}(A)}.
    \end{equation}
    Note that {\bf A1} implies, as in (\ref{Jmin}), the existence of
    a positive and finite
    constant $C_1$ ---that does not depend on $\mu$--- such that
    for all $\e>0$,
    \begin{equation}\label{eq:EZ}
        \E \left[ J_\e(\mu) \right] \ge C_1.
    \end{equation}

    Next, we will estimate the second moment
    of $J_\e(\mu)$. Let
    \begin{equation} \label{g}
    g_{\epsilon}(z):=\frac{1}{(2\e)^d} \1_{\tilde{B}(0,\e)}(z).
    \end{equation}
    Because
     \begin{equation}
        J_\e(\mu) = \int_{I} dt \, \int_{J} dx \, (g_{\epsilon}*\mu)
        (v(t\,,x)),
    \end{equation}
    Lemma \ref{prel}(1) and {\bf A2} together imply that there exists
    a finite and positive constant $C_2$ such that for all
    $\e>0$, 
\begin{align} \nonumber
        \E\left[ \left( J_\e(\mu) \right)^2\right]&= 
        \int_I dt \, \int_J dx \, \int_I ds \, \int_J dy \,
        \int_{\tilde{B}(z,\e)} dz_1  \int_{\tilde{B}(z,\e)} dz_2
        \\ \nonumber
        &\qquad \qquad\qquad \qquad \times p_{t,x;s,y}(z_1,z_2)\, (g_{\epsilon}*\mu)(z_1)\,(g_{\epsilon}*\mu)(z_2)\\
        &\le C_2 I_{\beta-6}(g_{\epsilon}*\mu).
\label{case2}
\end{align}
    By appealing to Theorem \ref{th:energybound:1}
    in Appendix \ref{apendA}, we see that  for all $\e>0$,
	\begin{equation}\begin{split}
        \E\left[ \left( J_\e(\mu) \right)^2\right]
       		&\le C_2 I_{\beta-6}(\mu)\\
			&\leq \frac{2C_2}{\textnormal{Cap}_{\beta-6}(A)}, 
	\end{split}\end{equation}
	by \eqref{eq:mu1}.
    The preceding, \eqref{eq:EZ}, and  the Paley--Zygmund inequality together imply that
    \begin{equation}
        \P\left\{ J_\e(\mu)>0 \right\}
        \ge \frac{C_1^2}{2C_2}\, \Cap_{\beta-6}(A).
    \end{equation}
    The left-hand side is bounded above by
    $\P\{ v(I \times J) \cap A^{(\e)}\neq\varnothing\}$,
    where $A^{(\e)}$ denotes the closed $\e$-enlargement
    of $A$. Let $\e\downarrow 0$ and appeal to the continuity
    of the trajectories of $v$ to find that for all $\mu \in \mathcal{P}(A)$,
    \begin{equation}
        \P\left\{ v(I\times J) \cap A \neq \varnothing\right\}
        \ge \frac{C_1^2}{2C_2}\, \Cap_{\beta-6}(A).
    \end{equation}

    \vskip 12pt

	{\em Case 3:}  $\beta-6=0$. We proceed as we
	did in Case 2, but use (\ref{case2}) with $\beta=6$
	and Theorem \ref{th:energybound:2} in the Appendix
	to obtain that  for all $\e>0$,
	\begin{equation}\begin{split}
        \E\left[ \left( J_\e(\mu) \right)^2\right]
        	&\le C_2 I_{0}(g_{\epsilon}*\mu)\\
        &\leq  c I_{0}(\mu)\\
        &\leq\frac{c}{\textnormal{Cap}_{0}(A)}. 
	\end{split}\end{equation}
	This proves part (1) of the theorem. 

     We prove (2) similarly. Without loss of generality
    we assume that $\Cap_{\beta-2}(A)>0$. This implies that $\beta-2 <d$ and $A \neq \varnothing$. Again, we need to consider three different cases.
\vskip 12pt

    {\em Case (i): $\beta-2<0$.} We proceed as we did in Case 1, but instead of
    $J_\e(z)$, we consider
    \begin{equation}
        \hat{J}_{\e,t}(z) := \frac{1}{(2\e)^d}
        \int_{J} dx \,  \1_{\tilde{B}(z,\e)}(v(t\,,x)),
    \end{equation}
for $t \in I$ fixed. We then use ${\bf A1}$ in order to obtain 
\begin{equation}
    \E\left[ \hat{J}_{\e,t}(z) \right] \geq C\, \vert J\vert >0.
\end{equation}
Note that, in this case, the constant $C$ depends on 
    $t$ only through $I$.
    We use ${\bf A2}$ to 
    bound the second moment of $\hat{J}_{\e,t}(z)$, that is,
\begin{equation}\begin{split}
        \E \left[ (\hat{J}_{\e,t}(z))^2 \right] &=
          \frac{1}{(2\e)^{2d}} \int_{J} dx \int_{J} dy 
          \int_{\tilde{B}(z,\e)} dz_1  \int_{\tilde{B}(z,\e)} dz_2 \,
          p_{t,x;s,y}(z_1,z_2) \\
        &\leq C  \int_0^{\vert J \vert} dv \,v^{-\beta/2},
\end{split}\end{equation}
    which is finite because $0<\beta<2$. The rest of the proof follows exactly as in Case 1. 
\vskip 12pt

    {\em Case (ii): $0<\beta-2<d$.} We choose
    $\mu\in\mathcal{P}(A)$ such that
    $I_{\beta-2}(\mu) \le 2/\Cap_{\beta-2}(A)$.
    We proceed as we did in Case 2, but instead of
    $J_\e(\mu)$, we consider
    \begin{equation}
        \hat{J}_{\e,t}(\mu) := \frac{1}{(2\e)^d}
        \int_{\R^d} \mu(dz) \, \int_{J} dx \, \1_{\tilde{B}(z,\e)}(v(t\,,x)),
    \end{equation}
    for $t \in I$ fixed. We then use ${\bf A1}$ 
in order to obtain 
\begin{equation}
   \E[\hat{J}_{\e,t}(\mu)] \geq C_1 >0.
\end{equation}
Finally, ${\bf A2}$ and Lemma \ref{prel}(2) with $\alpha=1$ and $I$ replaced by $J$ together imply that there exists a finite and positive constant $C$ such that for all $\e>0$,
    \begin{equation} 
        \E\left[ \left( \hat{J}_{\e,t}(\mu) \right)^2\right]
        \le C I_{\beta-2}(g_{\epsilon}*\mu).
    \end{equation}
     The remainder of the proof of (ii) follows exactly as we did for Case 2. 
\vskip 12pt

     {\em Case (iii): $\beta=2$.} We proceed as in (ii) and Case 3. This
     proves part (2) of the theorem.
\vskip 12pt

   We prove (3) by applying the same argument, but instead of
    $J_\e(\mu)$ and/or $\hat{J}_{\e,t}(\mu)$, consider
    \begin{equation}
        \bar{J}_{\e,x}(\mu) := \frac{1}{(2\e)^d}
        \int_{\R^d}\mu(dz) \, \int_{I} dt \, \1_{\tilde{B}(z,\e)}(v(t\,,x)),
    \end{equation}
    for $x \in J$ fixed, and use {\bf A1}, {\bf A2} 
    and Lemma \ref{prel}(2) with $\alpha=1/2$ to conclude.
\end{proof}

Theorem \ref{thm:cap:LB} is a result about hitting probabilities
of the random sets that are obtained by considering various
images of $v$. Next, we describe similar results for other,
related, random sets.
Define
\begin{equation}\begin{split}
    \mathcal{L}(z\,;v) &:= \left\{ (t\,,x)\in I \times
        J:\ v(t\,,x) = z\right\},\\
    \mathcal{T}(z\,;v) &= \left\{ t\in I:\ v(t\,,x)=z
        \text{ for some }x\in J \right\},\\
    \mathcal{X}(z\,;v) &= \left\{ x\in J:\ v(t\,,x)=z
        \text{ for some }t\in I\right\},\\
    \mathcal{L}_{x}(z\,;v) &:= \left\{ t\in I:\
        v(t\,,x)=z\right\},\\
    \mathcal{L}^{t}(z\,;v) &:= \left\{ x\in J:\
        v(t\,,x)=z\right\}.
\end{split}\end{equation}
We note that $\mathcal{L}(z\,;v)$ is the level set of $v$ at level $z$, $\mathcal{T}(z\,;v)$ (resp.~$\mathcal{X}(z\,;v)$) is the projection of $\mathcal{L}(z\,;v)$ onto $I$ (resp.~$J$), and $\mathcal{L}_{x}(z\,;v)$ (resp.~$\mathcal{L}^{t}(z\,;v)$) is the $x$-section (resp.~$t$-section) of $\mathcal{L}(z\,;v)$.

\begin{theorem}\label{thm:cap2:LB}
    Assume that {\bf A1} and {\bf A2} are met. Then, for all $R>0$,
    there exists a positive and finite constant $a=a(I,J,\beta,R,d)$ such that
    the following holds for all compact sets
    $E\subseteq I\times J$,
    $F\subseteq I$, and
    $G\subseteq J$, and for all $z\in B(0\,,R)$:
    \begin{enumerate}
        \item[\textnormal{(1)}]
            $\P\{ \mathcal{L}(z\,;v) \cap E
			\neq\varnothing \}\ge
			a \Cap_{\beta/2}^{{\bf \Delta}}
			(E)$;
        \item[\textnormal{(2)}]
            $\P\{ \mathcal{T}(z\,;v)\cap F\neq\varnothing\}\ge
			a \Cap_{(\beta-2)/4}(F)$;
        \item[\textnormal{(3)}]
            $\P\{\mathcal{X}(z\,;v)\cap G\neq\varnothing\}\ge
			a \Cap_{(\beta-4)/2} (G)$;
        \item[\textnormal{(4)}] 
             for all $x \in J$, $\P\{\mathcal{L}_x(z\,;v) \cap F\neq\varnothing\}\ge
			a\Cap_{\beta/4}(F)$;
        \item[\textnormal{(5)}] 
	        for all $t\in I$, $\P\{ \mathcal{L}^t(z\,;v) \cap G\neq\varnothing\}
			\ge a\Cap_{\beta/2}(G)$.
    \end{enumerate}
\end{theorem}

\begin{proof} 
    We begin by proving (1). Without loss of generality
    we assume that $\Cap_{\beta/2}^{{\bf \Delta}}(E)>0$.
    Choose $\mu\in\mathcal{P}(E)$ such that
    $I_{\beta/2}^{{\bf \Delta}}(\mu)\le 2/\Cap_{\beta/2}^{
    {\bf \Delta}}(E)$.
    For all $\delta>0$, define
    \begin{equation}
        Z_\delta(\mu) := \frac{1}{(2\delta)^d}
        \int_E \mu(dt\,dx) \, \1_{\tilde{B}(z,\delta)}(v(t\,,x)).
    \end{equation}
    Then, in accord with {\bf A1}, there exists a finite and positive
    constant $C_1$ such that for all $\mu\in\mathcal{P}(E)$ and $\delta>0$,
    \begin{equation}\label{eq:bd1}
        \E[Z_\delta(\mu)] \ge C_1.
    \end{equation}
    On the other hand, {\bf A2} guarantees the
    existence of a finite and positive constant $C_2$
    such that for all $\mu\in\mathcal{P}(E)$ and $\delta>0$,
    \begin{equation}\label{eq:bd2}\begin{split}
        \E\left[ \left( Z_\delta(\mu) \right)^2\right] 
			&= \frac{1}{(2 \delta)^{2d}}
	        \int_E \mu(dt\, dx) \int_E \mu(ds \, dy) 
	        \int_{\tilde{B}(z\,, \delta)} dz_1  \int_{\tilde{B}(z\,, \delta)} dz_2 \
	        p_{t,x;s,y}(z_1\,,z_2) \, \\
		&\le C_2\int_E \int_E \frac{\mu(dt\, dx)\, \mu(ds\, dy)}{
	        \left[ {\bf \Delta}\left((t\,,x)\,; (s\,,y)
	        \right) \right]^{\beta/2}}\\
		&\le \frac{2C_2}{\Cap_{\beta/2}^{{\bf \Delta}}(E)}.
    \end{split}\end{equation}
    Equations \eqref{eq:bd1} and \eqref{eq:bd2},
    together with the Paley--Zygmund inequality, imply that
    \begin{equation}
        \P\left\{ Z_\delta(\mu) >0 \right\} \ge
        \frac{C_1^2}{2C_2}\, \Cap_{\beta/2}^{{\bf \Delta}}
        (E).
    \end{equation}
    The left-hand side is clearly bounded above by
    \begin{equation}
	    \P\left(\bigcup_{z_1\in B(z,\delta)}
	    (\mathcal{L}(z_1\,;v)\cap E) \neq\varnothing\right).
	\end{equation}
    Let $\delta\downarrow 0$ to finish the proof of (1).

    In order to prove (2), define, for all $\mu\in\mathcal{P}(F)$,
    $\delta>0$ and $z \in B(0,R)$,
    \begin{equation}
        Z_\delta(\mu) = \frac{1}{(2\delta)^d}
        \int_F \mu(dt) \, \int_J dx\, \1_{\tilde{B}(z,\delta)}(
        v(t\,,x)).
    \end{equation}
    By \textbf{A1}, we can find a constant $C$ ---
    depending only on $(I\,,J\,,R\,,d)$ --- such that
    \begin{equation}
	    \inf_{\delta>0}\ 
	    \inf_{\mu\in\mathcal{P}(F)}
	    \E \left[ Z_\delta(\mu) \right] \ge C.
    \end{equation}
    On the other hand, let $g_{\delta}$ be as defined in
    (\ref{g}) with $\e$ replaced by $\delta$. By \textbf{A2},
    there exists $\tilde{C}$---depending only on
    $(I\,,J\,,\beta\,,R\,,d)$---such that for all $\delta>0$
    and $\mu\in\mathcal{P}(F)$,
\begin{equation}\begin{split}
	\E\left[ \left( Z_\delta(\mu) \right)^2\right] &=
		\int_F \mu(dt) \, \int_J dx\, \int_F \mu(ds) \,
		\int_J dy\, \int_\R dz_1 \int_\R dz_2\, \\
	&\hskip2in\times g_{\delta}(z_1 - z)
		g_{\delta}(z_2 - z)\, p_{t,x;s,y}(z_1\,,z_2).
\end{split}\end{equation}
Since
\begin{equation}
        \int_J dx \, \int_J dy \,
        \frac{1}{\left[ {\bf \Delta}
        \left( (t\,,x)\,;(s\,,y)\right)\right]^{
        \beta/2}} \le 2 \vert J \vert \Psi_{\vert J \vert,\beta/2}
        \left( |t-s|^{\beta/4}\right),
\end{equation}
where $\Psi_{a,\nu}(\rho)$ is defined in (\ref{rdpsi}), we see that
\begin{equation}\begin{split}
    \E\left[ \left( Z_\delta(\mu) \right)^2\right]
    &\le C \int_F \mu(dt)\, \int_F \mu(ds) \, \int_\R dz_1 \int_\R dz_2\\
    &\hskip1.7in
    	\times g_{\delta}(z_1 - z) g_{\delta}(z_2 - z)\,
		\Psi_{\vert J \vert,\beta/2}(\vert t-s\vert^{\beta/4}).
\end{split}\end{equation}
Since the two $dz_i$-integrals are equal to $1$, 
Lemma \ref{lem:Psi} implies that there exists a
constant $\bar{C}$ such that for all $\mu\in\mathcal{P}(F)$ and $\delta>0$, 
\begin{equation}\begin{split}
	\E\left[ \left( Z_\delta(\mu) \right)^2\right]
		&\le \bar{C}\int_F \mu(dt) \, \int_F \mu(ds) \,
		{\rm K}_{1-(2/\beta)}(|t-s|^{\beta/4})\\
	&= \bar{C} I_{(\beta-2)/4}(\mu).
\end{split}\end{equation}
An application of the Paley--Zygmund inequality implies statement (2) of the theorem.

    In order to prove (3), we consider instead $\mu\in\mathcal{P}(G)$ and
    \begin{equation}
        \bar{Z}_\delta(\mu) = \frac{1}{(2\delta)^d}
        \int_G  \mu(dx) \, \int_I dt \ \1_{\tilde{B}(z,\delta)}(
        v(t\,,x)).
    \end{equation}
    Thanks to {\bf A1}, $\E[\bar{Z}_\delta(\mu)]$ is bounded below,
    uniformly for all $\delta >0$ and
    $\mu\in\mathcal{P}(G)$. Also, as above, {\bf A2} implies that
    there exists a positive and finite constant $C$ such that
    $\E[(\bar{Z}_\delta(\mu))^2]\le C I_{(\beta-4)/2}(\mu)$
    for all $\delta >0 $ and $\mu\in\mathcal{P}(G)$.
    Indeed, this is a consequence of Lemma \ref{lem:Psi}
    and the fact that
    \begin{equation}
        \int_I dt \, \int_I ds \, 
        \frac{1}{\left[ {\bf \Delta}
        \left( (t\,,x)\,;(s\,,y)\right)\right]^{
        \beta/2}} \le 2 \vert I \vert \Psi_{\vert I \vert,\beta/4}
        \left( |x-y|^{\beta/2}\right).
    \end{equation}
Therefore, statement (3) now follows from 
the two moment bounds and the Paley--Zygmund inequality.

    For (4), we consider instead $z \in B(0,R)$, $ x \in J$, $\mu\in\mathcal{P}(F)$ and set
\begin{equation}
        Z'_\delta(\mu) = \frac{1}{(2\delta)^d}
        \int_F \mu(dt)\, \1_{\tilde{B}(z,\delta)}(v(t\,,x)) .
\end{equation}
As was the case in (1), (2), and (3), $\E[Z'_\delta(\mu)]$ is bounded below, uniformly for all $\delta >0$, $\mu\in\mathcal{P}(F)$ and $x \in J$. In addition, there exists a positive and finite constant $C$ such that 
\begin{equation}
   \E[(Z'_\delta(\mu))^2] \le \int_F \mu(dt)\, \int_F \mu(ds)\, \int_\R dz_1 \int_\R dz_2\, g_{\delta}(z_1 - z) g_{\delta}(z_2 - z)\, p_{t,x;s,x}(z_1\,,z_2).
\end{equation}
Since $p_{t,x;s,x}(z_1\,,z_2) \le \vert t-s\vert^{-\beta/4}$, and the two $dz_i$-integrals are equal to $1$, we see that
\begin{equation}
   \E[(Z'_\delta(\mu))^2] \le C I_{\beta/4}(\mu),
\end{equation}
for all $\delta >0$, $\mu\in\mathcal{P}(F)$ and $x \in J$.
Therefore, statement (4) follows from 
the two moment bounds and the Paley--Zygmund inequality.

    Finally, in order to prove (5), we consider instead
    $\mu\in\mathcal{P}(G)$ and
    \begin{equation}
        Z''_\delta(\mu) = \frac{1}{(2\delta)^d}
        \int_G \1_{\tilde{B}(z,\delta)}(v(t\,,x))\, \mu(dx).
    \end{equation}
    Once again by {\bf A1}, $\E[Z''_\delta(\mu)]$ is bounded below,
    uniformly for all $\delta>0$ and $\mu\in\mathcal{P}(F)$.
    And by {\bf A2},
    $\E[(Z''_\delta(\mu))^2]\le C I_{\beta/2}(\mu)$, where
    $C\in\,]0\,,\infty[$ does not depend on $(\delta\,,\mu)$. From the two moment bounds,
(5) follows, whence the theorem.
\end{proof}

\begin{remark} (a) Hypothesis {\bf A1} is convenient since, together with {\bf A2}, it leads to all the conclusions of Theorem \ref{thm:cap:LB} and \ref{thm:cap2:LB}. If one is only interested in certain of these conclusions, then weaker assumptions than {\bf A1} are possible, analogous to Hypothesis H1 of \cite{Dalang:04}. For instance, Theorem \ref{thm:cap:LB}(1) can be obtained if {\bf A1} is replaced by:
\vskip 12pt

\noindent{\bf A1'}. For all $M>0$, there exists a positive and finite constant $C = C(I,J,M,d)$ such that for all $z\in[-M\,,M]^d$,
\begin{equation}
           \int_I dt \int_J dx\,  p_{t,x}(z) \ge C.
\end{equation}
Indeed, this assumption would be used to get the lower bound in (\ref{Jmin}) and (\ref{eq:EZ}).

   In the same way, Theorem \ref{thm:cap:LB}(2) can be obtained if {\bf A1} is replaced by:
\vskip 12pt

\noindent{\bf A1$^t$}. For all $M>0$, there exists a positive and finite constant $C = C(t,J,M,d)$ such that for all $z\in[-M\,,M]^d$,
\begin{equation}
           \int_J dx\, p_{t,x}(z) \ge C.
\end{equation}

   Similar considerations apply to Theorem \ref{thm:cap:LB}(3), which can be obtained if {\bf A1} is replaced by:
\vskip 12pt

\noindent{\bf A1$_x$}. For all $M>0$, there exists a positive and finite constant $C = C(x,I,M,d)$ such that for all $z\in[-M\,,M]^d$,
\begin{equation}
           \int_I dt\, p_{t,x}(z) \ge C.
\end{equation}

   (b) It is also possible to weaken Hypothesis {\bf A2}. For instance, Theorems \ref{thm:cap:LB}(2) and \ref{thm:cap2:LB}(5) can be proved if {\bf A2} is replaced by:
\vskip 12pt

\noindent{\bf A2$^t$}. There exists $\beta>0$ such that for all $M>0$, there exists $c = c(t,I,J,\beta,M,d)>0$ such that for all $x,y\in J$ with $x \neq y$, and for every $z_1,z_2\in[-M\,,M]^d$,
        \begin{equation}
            p_{t,x;t,y}(z_1\,,z_2) \le
            \frac{c}{ \vert x - y \vert^{\beta/2}}\exp\left(
            - \frac{\| z_1-z_2\|^2}{c
            \vert x - y \vert}\right).
        \end{equation}
Similar considerations also apply to Theorem \ref{thm:cap:LB}(3).
\end{remark}

\section{Upper Bounds on Hitting Probabilities}\label{sec3}

The results of this section complement those of the
preceding by establishing upper bounds for various
hitting probabilities.

Consider two compact
nonrandom intervals $I\subset [0\,,T]$ 
and $J\subset [0\,,1]$, and suppose
${v}=\{{v}(t\,,x)\}_{(t,x)\in I\times J}$
is an $\R^d$-valued random field.
For all positive integers $n$, set $t_k^n:=k 2^{-4n}$, $x_{\ell}^n:=\ell 2^{-2n}$,
and
\begin{equation} \label{covering}
	I^n_k = [t_k^n,t_{k+1}^n],\qquad J^n_{\ell} = [x_{\ell}^n,x_{\ell+1}^n],\qquad R^n_{k,\ell} = I^n_k \times J^n_{\ell}.
\end{equation}

\begin{theorem}\label{thm:meas:UB} 
    Fix $\beta>0$ and $M>0$.
Suppose that there exists $c >0$ such that for
	all $z \in [-M,M]^d$, $\epsilon>0$, large $n$ and $R_{k,\ell}^n \subseteq I\times J$, 
	\begin{equation}\label{cond:i}
		\P \{ v(R_{k,\ell}^n) \cap B(z\,,\e) \neq \varnothing\} \le c\, \e^{\beta}.
	\end{equation}
    Then there exists a positive and finite constant $a$ such that for all Borel sets $A \subset [-M,M]^d$:
\begin{enumerate}
   \item[\textnormal{(1)}] $\P\{ {v}(I\times J)\cap A
   \neq\varnothing \} \le a \mathcal{H}_{\beta-6}(A)$;
  \item[\textnormal{(2)}] 
  for every $t\in I$,  $\P\{ {v}(\{t\}\times J) \cap A \neq\varnothing \}
            \le a \mathcal{H}_{\beta-2}(A)$;
  \item[\textnormal{(3)}] 
  for every $x\in J$,  $\P\{ {v}(I\times \{x\}) \cap A \neq\varnothing \}
            \le a \mathcal{H}_{\beta-4}(A)$.
\end{enumerate}
\end{theorem}

\begin{proof} We begin by proving (1). When $\beta - 6 <0$, there is nothing to prove, so we assume that $\beta - 6 \geq 0$.
	Fix $\e \in\, ]0\,,1[$ and $n \in \mathbb{N}$ such that 
	$2^{-n-1} <\epsilon \leq 2^{-n}$, and write
    \begin{equation}
	    \P\left\{ v\left( I\times J\right)
		\cap B(z\,,\e) \neq\varnothing \right\}
		\le \mathop{\sum\sum}_{%
		\substack{(k,\ell):\\
	   R_{k,\ell}^n\cap( I\times J) \neq \varnothing}}
		\P \{
		v(R_{k,\ell}^n) \cap B(z\,,\e) \neq \varnothing\}.
    \end{equation}
    The number of pairs $(k\,,\ell)$ involved in the two 
    sums is at most $2^{6n}$. Because $2^{-n-1}< \epsilon$,
    the condition \eqref{cond:i} implies that
    for all large $n$ and all $z \in A$,
    \begin{equation}\label{eq:hit:ball:UB}\begin{split}
        \P\left\{ {v}\left( I\times J\right)
        	\cap B(z\,,\e) \neq\varnothing \right\}
        	&\le \tilde{C} 2^{-n(\beta-6)}\\
			&\le C \e^{\beta-6}.
    \end{split}\end{equation}
    Note that $C$ does not depend on $(n\,,\e)$.
    Therefore, (\ref{eq:hit:ball:UB}) is valid for all
    $\epsilon \in\, ]0\,,1[$.

    Now we use a \emph{covering argument}:
    Choose $\epsilon \in\,]0\,,1[$ and let $\{B_i\}_{i=1}^\infty$ 
    be a sequence of open balls in $\R^d$ with
    respective radii $r_i\in\,]0\,,\e]$ such that
    \begin{equation}\label{eq:HHHH}
	    A\subseteq \bigcup_{i=1}^\infty B_i
	    \quad\text{and}\quad
	    \sum_{i=1}^\infty (2r_i)^{\beta-6} \le
	    \mathcal{H}_{\beta-6} (A)+\e.
	\end{equation}
Because $\P\{ {v}( I\times J)
	\cap A \neq\varnothing\}$ is at most 
	$\sum_{i=1}^\infty \P\{
	{v}(I\times J)\cap B_i\neq
	\varnothing\}$, \eqref{eq:hit:ball:UB} 
	and \eqref{eq:HHHH} together imply that
    \begin{equation}\begin{split}
			\P\left\{ {v}\left( I\times J\right)
				\cap A \neq\varnothing \right\}
				&\le C \sum_{i=1}^\infty r_i^{\beta-6}\\
        &\le \bar{C} ( \mathcal{H}_{\beta-6}
        	(A)+\e).
	\end{split}\end{equation}
    Let $\e\to 0^+$ to deduce (1).

    In order to prove (2), we can assume that
    $\beta - 2 \geq 0$ and we fix $\epsilon \in\,]0\,,1[$.
    We can find integers $n$ and $k$ such that $2^{-n-1}<\epsilon \leq 2^{-n}$
    and $t \in I^n_k$. Then, by \eqref{cond:i},
    \begin{equation}\begin{split}
        \P\left\{ {v}\left(\{t\}\times J\right)
				\cap B(z\,,\e) \neq\varnothing \right\} 
				&\le  \sum_{\ell: J^n_{\ell} \cap J\neq \varnothing}
				\P \{ v(I_{k}^n \times J^n_{\ell}) \cap B(z\,,\e) \neq \varnothing\} \\
			&\le C 2^{-n\beta} 2^{2n}\\
			&\leq \tilde{C} \e^{\beta-2}.
    \end{split}\end{equation}
    Now use a covering argument, as we did to prove (1), in order
    to verify (2).

    The proof of (3) follows along similar lines, and is left to the reader.
\end{proof}

\begin{theorem}\label{thm:cap2:LBB}
	Fix $\beta >0$ and $M>0$.
    If the assumptions of Theorem \ref{thm:meas:UB} are met, then
    there exists $a\in\,]0\,,\infty[$ such that
    the following holds for all $z\in [-M,M]^d$ and all compact sets
    $E\subseteq I \times J$, $F\subseteq I$, and $G\subseteq J$:
    \begin{enumerate}
        \item[\textnormal{(1)}] $\P\{ \mathcal{L}(z\,;v) \cap E
			\neq\varnothing \}\le
			a \mathcal{H}^{{\bf \Delta}}_{\beta/2}(E)$;
        \item[\textnormal{(2)}]
            $\P\{ \mathcal{T}(z\,;v)\cap F\neq\varnothing\}\le
			a \mathcal{H}_{(\beta-2)/4}(F)$;
        \item[\textnormal{(3)}]
            $\P\{\mathcal{X}(z\,;v)\cap G\neq\varnothing\}\le
			a  \mathcal{H}_{(\beta-4)/2}(G)$;
        \item[\textnormal{(4)}]
           for all $x \in J$,  $\P\{\mathcal{L}_x(z\,;v) \cap F\neq\varnothing\}\le
			a \mathcal{H}_{\beta/4}(F)$;
        \item[\textnormal{(5)}] 
          for all $t\in I$, $\P\{ \mathcal{L}^t(z\,;v) \cap G\neq\varnothing\}
			\le a \mathcal{H}_{\beta/2}(G)$.
    \end{enumerate}
\end{theorem}

\begin{proof} 
	 Let $z \in [-M,M]^d$. Fix $r \in\,]0\,,1[$, $t_0 \in I$
	 and $x_0 \in J$. We can find integers $n$, $\ell$ and $k$
	 such that $2^{-2n-2}<r \leq 2^{-2n-1}$, $t_0 \in I_{k}^n$,
	 $x_0 \in J_{\ell}^n$. Then condition (\ref{cond:i})
	 implies that for $n$ large,
	\begin{equation} \label{aux1}\begin{split}
		\P\left\{ \inf_{\substack{t_0 \leq t \leq t_0+r^{1/2}\\
			x_0 \le x\le x_0+r}} |v(t\,,x)-z|\le r^{1/2} \right\} &\le
			\sum_{i=k}^{k+1} \sum_{j=\ell}^{\ell+1}
			\P \{ v(R_{i,j}^n) \cap B(z\,,r^{1/2}) \neq \varnothing\} \\
			&\le 
			C r^{\beta/2}.
	\end{split}\end{equation}
Note that $C$ does not depend on $(n\,,r\,,t_0\,,x_0)$. 
	 
	 Now we use a \emph{covering argument}:
    Choose $r \in\,]0\,,1[$ and let $\{E_i\}_{i=1}^\infty$ 
    denote a sequence of open ${\bf  \Delta}$-balls in $I \times J$ with
    respective radii $r_i\in\,]0\,,r]$ such that
    \begin{equation}
	    E\subseteq \bigcup_{i=1}^\infty E_i
	    \quad\text{and}\quad
	    \sum_{i=1}^\infty (2r_i)^{\beta/2} \le
	    \mathcal{H}^{{\bf \Delta}}_{\beta/2} (E)+r.
	\end{equation}
	Then
    \begin{equation}\begin{split}
        \P\left\{  \mathcal{L}(z\,;v)\cap E\neq\varnothing \right\}
            &= \P\left\{ \inf_{(t,x)\in E}
            |v(t\,,x)-z|=0 \right\}\\
        &\le \sum_{i=1}^\infty \P\left\{ \inf_{(t,x)\in E_i}
            |v(t\,,x)-z|\le r_i^{1/2}\right\}\\
        &\le C \sum_{i=1}^\infty r_i^{\beta/2}\\
			&\le \tilde{C} (\mathcal{H}^{{\bf \Delta}}_{\beta/2} (E)+r).
    \end{split}\end{equation}
    Let $r\to 0^+$ to deduce (1).

	To prove (2), fix $r  \in\,]0\,,1[$ and $t_0 \in I$.
	There exist integers $n$ and $k$ such that
	$2^{-4n-2}<r \leq 2^{-4n-1}$ and $t_0 \in I_k^n$.
	Condition (\ref{cond:i}) implies that for $n$ large,
    \begin{equation} \label{aux2}\begin{split}
			\P\left\{ \inf_{t_0 \leq t \leq t_0+r} \ 
				\inf_{x \in J} |v(t\,,x)-z|\le r^{1/4} \right\} &\le
        	\sum_{i=k}^{k+1}\ 
				\mathop{\sum}_{\ell:
				J_{\ell}^n\cap J \neq \varnothing}
				\P \{ v(R_{i,\ell}^n) \cap B(z\,,r^{1/4}) \neq \varnothing\}\\
        &\le \tilde{C} 2^{-n\beta} 2^{2n}\\
        &\leq  C r^{(\beta-2)/4},
    \end{split}\end{equation}
since $r>2^{-4n-2}$. Note that $C$ does not depend on $(n\,,r\,, t_0)$. 

    Choose $r \in\,]0\,,1[$ and let $\{F_i\}_{i=1}^\infty$ 
    denote a sequence of open balls in $I$ with
    respective radii $r_i\in\,]0\,,r]$ such that
    \begin{equation}
	    F\subseteq \bigcup_{i=1}^\infty F_i
	    \quad\text{and}\quad
	    \sum_{i=1}^\infty (2r_i)^{(\beta-2)/4} \le
	    \mathcal{H}_{(\beta-2)/4} (F)+r.
	\end{equation}
	Then
    \begin{equation}\begin{split}
        \P\left\{  \mathcal{T}(z\,;v)\cap F\neq\varnothing \right\}
            &= \P\left\{ \inf_{t\in F}\ \inf_{x\in J}
            |v(t\,,x)-z|=0 \right\}\\
        &\le \sum_{i=1}^\infty \P\left\{ \inf_{t\in F_i}\ \inf_{x\in J}
            |v(t\,,x)-z|\le r_i^{1/4}\right\}\\
        &\le C \sum_{i=1}^\infty r_i^{(\beta-2)/4}\\
        & \le \tilde{C}
        	(\mathcal{H}_{(\beta-2)/4} (E)+r).
    \end{split}\end{equation}
     Let $r\to 0^+$ to deduce (2).

	The proof of (3) follows along similar lines, and is left to the reader.

	We now prove (4). Fix $x \in J$, $r \in\,]0\,,1[$ and $t_0 \in I$.
	There exist integers $n$, $k$ and $\ell$ such that
	$2^{-4n-2}<r \leq 2^{-4n-1}$, $t_0 \in I_k^n$ and $x \in J_{\ell}^n$. Condition (\ref{cond:i}) implies that for $n$ large,
	\begin{equation} \label{aux3}\begin{split}
        \P\left\{ \inf_{t_0 \leq t \leq t_0+r} 
	        |v(t\,,x)-z|\le r^{1/4} \right\} &\le \sum_{i=k}^{k+1}
            \P \{ v(R_{i,\ell}^n) \cap B(z\,,r^{1/4}) \neq \varnothing\} \\
           & \le 
              C r^{\beta/4},
    \end{split}\end{equation}
Note that $C$ does not depend on $(n\,,r\,,x, t_0)$. 

    Choose $r \in\,]0\,,1[$ and let $\{F_i\}_{i=1}^\infty$ 
    denote a sequence of open balls in $I$ with
    respective radii $r_i\in\,]0\,,r]$ such that
    \begin{equation}
	    F\subseteq \bigcup_{i=1}^\infty F_i
	    \quad\text{and}\quad
	    \sum_{i=1}^\infty (2r_i)^{\beta/4} \le
	    \mathcal{H}_{\beta/4} (F)+r.
	\end{equation}
Then
    \begin{equation}\begin{split}
        \P\left\{  \mathcal{L}_x(z\,;v)\cap G \neq\varnothing \right\}
            &= \P\left\{ \inf_{t\in F}
            |v(t\,,x)-z|=0 \right\}\\
        &\le \sum_{i=1}^\infty \P\left\{ \inf_{t\in F_i}
            |v(t\,,x)-z|\le r_i^{1/4}\right\}\\
        &\le C \sum_{i=1}^\infty r_i^{\beta/4}\\
        &\le \tilde{C} (\mathcal{H}_{\beta/4} (E)+r).
    \end{split}\end{equation}
    Let $r\to 0^+$ to deduce (4).

The proof of (5) follows along similar lines, and is left to the reader.
\end{proof}

The results of this section all assume Condition (\ref{cond:i}). The following provides a useful sufficient condition for (\ref{cond:i}) to hold. This conditions is used for instance in \cite{Dalang-K-N:07}.

\begin{theorem} \label{eta}
	Fix $M>0$. Assume that the $\R^d$-valued
	random field $v$ satisfies the following two conditions:
	\begin{itemize}
		\item[\textnormal{(i)}] For any $(t\,,x) \in I \times J$,
			the random vector $v(t\,,x)$ has a density $p_{t,x}(z)$ which
			is is uniformly bounded over $z \in [-M,M]^d$ and $(t\,,x) \in I\times J$.

		\item[\textnormal{(ii)}] For all $p>1$, there exists a 
			constant C depending on $p, I, J$ such that for
			any $(t\,,x), (s\,,y) \in I \times J$, 
			\begin{equation} \label{Hi}
				\E [\vert v(t\,,x)-v(s\,,y) \vert^p] \leq C \left[ {\bf \Delta}
				\left( (t\,,x)\,;(s\,,y)\right) \right]^{p/2}.
			\end{equation}
	\end{itemize}
Then for any $\beta \in\,]0\,,d[$, Condition 
    \textnormal{(\ref{cond:i})} is satisfied  and therefore, so are the upper bounds on hitting probabilities in Theorems \ref{thm:meas:UB} and \ref{thm:cap2:LBB} for such $\beta$.
    \end{theorem}

 \begin{proof}
Fix $z\in [-M,M]^d$. For $n \in \mathbb{N}$ and $\varepsilon \in\,]0\,,1[$, set 
\begin{equation}\begin{split}
    Y^n_{k,\ell} &:= \left| {v}(t_k^n\,,x^n_{\ell}) - z\right|, \\
    Z^n_{k,\ell} &:= \sup_{(t,x)\in B_{{\bf \Delta}}
        ((t_k^n,x_{\ell}^n),\varepsilon^2)} \left| {v}(t\,,x) - {v}(t^n_k\,,x^n_{\ell})
        \right|.
\end{split}\end{equation}
Fix $\beta \in\,]0\,,d[$. We are going to start by showing that
\begin{equation}\label{tt}
   \P\left\{Z^n_{k,\ell} \geq \frac{1}{2} Y^n_{k,\ell}\right\}\le  \tilde{c} \varepsilon^{\beta}.
\end{equation}
Indeed, observe that
\begin{equation}
\P\left\{Z^n_{k,\ell} \geq \frac{1}{2} Y^n_{k,\ell}\right\} \leq \P\left\{Y^n_{k,\ell} \leq \varepsilon^{\beta/d}\right\} + \P\left\{Z^n_{k,\ell} \geq \frac{1}{2} \varepsilon^{\beta/d}\right\}.
\end{equation}
By hypothesis (i), the first term on the right-hand side is
bounded by $c \varepsilon^{\beta}$. By Markov's inequality,
\begin{equation}
   \P\left\{Z^n_{k,\ell} \geq \frac{1}{2}
   \varepsilon^{\beta/d}\right\} \le
   \left(\frac{1}{2} \varepsilon^{\beta/d} \right)^{-p}
   \E[\vert Z^n_{k,\ell} \vert^p].
\end{equation}
Let $p>6$ and $q=\frac{p}{2} - 3$. Then $q >0$ and
$\frac{q}{p} = \frac{1}{2} - \frac{3}{p} >0$. Since
$\frac{\beta}{2d} < \frac{1}{2}$, we can choose $p$
large enough that $\frac{1}{2} - \frac{3}{p} > \frac{\beta}{2d}$.

Fix $\alpha \in \,]\frac{\beta}{2d}\,,\frac{q}{p}[$.
By hypothesis (ii) and Corollary \ref{AKCT},
\begin{equation}
     \E\left( \vert Z^n_{k,\ell} \vert^p\right) 
     \le (\varepsilon^2)^{\alpha p},
\end{equation}
and hence,
\begin{align}
   \P\left\{Z^n_{k,\ell} \geq \frac{1}{2}
   	Y^n_{k,\ell}\right\} &\leq c \, \varepsilon^{\beta}
		+c\, \varepsilon^{2 \alpha p -\beta p/d} \\
	&\leq c \, \varepsilon^{\beta} (1+c
		\varepsilon^{p(2 \alpha-\beta/d)-\beta}).
\end{align}
Since $2 \alpha-\beta/d >0$, it follows that $p(2 \alpha-\beta/d)-\beta >0$
for all sufficiently large $p$.
This proves (\ref{tt}). 

   Now, let $\varepsilon \in\,]0\,,1[$ and $n \in \mathbb{N}$
   be such that $2^{-n-1}<
\varepsilon \leq 2^{-n}$. According to (\ref{tt}),
    \begin{equation}\begin{split}
        \P\left\{ v\left( R_{k,\ell}^n \right)
            \cap B(z\,,\varepsilon) \neq\varnothing \right\}
        &\le
            \P\left\{ Y_{k,\ell}^n \le 2^{-n}+Z^n_{k,\ell}  \right\}\\
        &\le \P\left\{ Z^n_{k,\ell} \ge \frac12 Y^n_{k,\ell} \right\}
            + \P\left\{ Y^n_{k,\ell} \le 2^{1-n} \right\} \\
        &\le  c 2^{-n\beta} + c 2^{(1-n)d}.
    \end{split}\end{equation}
Therefore, for all large $n$ and all $z \in [-M,M]^d$,
    \begin{equation}
        \P\left\{ v \left( R_{k,l}^n\right)
        \cap B(z\,,\varepsilon) \neq\varnothing \right\}
        \le C 2^{-n\beta} \le \tilde{C} \varepsilon^{\beta},
    \end{equation}
    since $2^{-n-1}< \varepsilon$. This proves (\ref{cond:i}) and whence the theorem.
 \end{proof}

\section{The Gaussian case}\label{sec4}

We consider the s.p.d.e.\ \eqref{eq:spde} in the drift-free case
($b_i\equiv 0$), and write it in vector notation as
\begin{equation} \label{e:p}
    \frac{\partial{u}}{\partial t} = \frac{
    \partial^2{u}}{\partial x^2} +
    {\s}\dot{{W}}.
\end{equation}
The solution is the $d$-dimensional Gaussian random field
$\{u(t\,,x)\}_{t\in[0,T],x\in[0,1]}$ defined by
\begin{equation} \label{rf}
    u(t\,,x) = \int_0^t \int_0^1 G_{t-s}(x\,,y)
    \,\sigma \, W(ds dy),\qquad
    0\le t\le T,\ 0\le x\le 1.
\end{equation}

The main objective of this section is to show that for
$t_0>0$, the conclusions of Theorems \ref{thm:cap:LB},
\ref{thm:cap2:LB}, \ref{thm:meas:UB}, and \ref{thm:cap2:LBB}
are satisfied for $(u(t\,,x))$ with $\beta = d$, $I = [t_0\,,T]$,
and $J=[0\,,1]$. We point out that it would be much simpler to establish this for $\beta < d$: see the comment just before Proposition \ref{haus}. We begin with the following.

\begin{proposition} \label{capa}
Fix $t_0 >0$. 	Then the solution to
	\textnormal{(\ref{e:p})} satisfies 
	\textbf{A1} and \textbf{A2}
	with $\beta=d$, $I=[t_0\,,T]$ and $J=[0\,,1]$.
\end{proposition}

\begin{proof}
	It suffices to prove that Hypotheses {\bf A1} and {\bf A2} are
	satisfied for the random field (\ref{rf}). We are going to reduce the problem to the case where $\sigma$ is the $d\times d$ identity matrix by a change of variables.
		Because $\bm{\s}$ is invertible,
\begin{equation*}
	\frac{\partial (\s^{-1} u)}{\partial t}
	= \frac{\partial^2(\s^{-1} u)}{\partial x^2}
	+ \dot{W}.
\end{equation*}
Define $v:= \s^{-1}u$
to find that $v$ solves the following
\emph{uncoupled} system of s.p.d.e.'s:
\begin{equation} \label{v}
	\frac{\partial v}{\partial t} = \frac{\partial^2 v}{
	\partial x^2} + \dot{W}.
\end{equation}
We will prove that Hypotheses {\bf A1} and {\bf A2} hold for the solution of (\ref{v}).
Therefore, they also hold for $u=\sigma v$.
Note that $v=(v_1,\dots,v_d)$, where $v_1,\ldots,v_d$ are i.i.d.~real-valued processes.

\vskip 12pt
	\noindent{\it Verification of {\bf A1}.} 
Fix $I=[t_0\,,T]$, $J=[0\,,1]$ and $M>0$,
and let $z \in [-M\,,M]^d$.
		Then, for all $(t\,,x) \in I\times J$, the 
		probability density function of $v(t\,,x)$ is given by
		\begin{equation}
			p_{t,x}(z)=
			\frac{1}{(2 \pi \sigma_{t,x}^2)^{d/2}} \exp
			\left(- \frac{\Vert z \Vert^2}{2
			\sigma_{t,x}^2} \right),
		\end{equation}
		where
		\begin{equation}
			\sigma_{t,x}^2 := \text{Var}\
			v_i(t\,,x) =\int_0^t dr \int_0^1 dv \,
			(G_{t-r}(x\,,v))^2.
		\end{equation}
Since $(t\,,x) \mapsto \sigma_{t,x}^2$ is a continuous function, it achieves its minimum $\rho_1 >0$ and its maximum $\rho_2 < \infty$ over $I\times J$.		Thus,
		\begin{equation}
			  p_{t,x}(z) 
			\geq \frac{1}{(2 \pi \rho_2)^{d/2}} 
			\exp\left( -\frac{M^2 d}{2 \rho_1} \right).
		\end{equation}
		This proves {\bf A1}.

	\vskip 12pt

\noindent{\it Verification of {\bf A2}.} 
	We follow the proof of Theorem 2.1 of \textsc{Dalang and Nualart} 
	\cite{Dalang:04}. The joint probability density
	function $ p^i_{t,x;s,y}(\cdot\,,\cdot)$ of
	$(v_i(t\,,x)\,, v_i(s\,,y))$---for any two distinct space-time points
	$(t\,,x)$ and $(s\,,y)$---does not depend on $i$ and can be written as
	\begin{equation}
			p^i_{t,x;s,y}(z_1\,,z_2)= p^i_{t,x\vert s,y}(z_1\mid z_2) p^i_{s,y}(z_2),
	\end{equation}
	where $z_1\,,z_2 \in \R$, $p^i_{t,x\vert s,y}(\ \cdot \mid z_2)$ denotes
	the conditional probability density function of $v_i(t\,,x)$ 
	given $v_i(s\,,y) = z_2$ and $p^i_{s,y}(\cdot)$ denotes the marginal
	density of $v_i(s\,,y)$. By linear regression,
	\begin{equation}
		p^i_{t,x\vert s,y}(z_1\mid z_2)= \frac{1}{ \tau \sqrt{ 2 \pi}  }
		\exp \left(- \frac{\vert z_1-m z_2 \vert^2}{2 \tau^2} \right) ,
	\end{equation}
	where
	\begin{equation}\begin{split}
		& \tau^2:=\tau^2_{t,x;s,y}=\sigma_{t,x}^{2} 
			\left(1-\rho^2_{t,x;s,y}\right), \qquad \rho_{t,x;s,y}=
			\frac{\sigma_{t,x;s,y}}{\sigma_{t,x}\sigma_{s,y}} \\
		& m:=m_{t,x;s,y}=\frac{\sigma_{t,x;s,y}}{\sigma_{s,y}^{2}}, \qquad
			\sigma_{t,x;s,y}= \text{Cov}\left(v_i(t\,,x)\,, v_i(s\,,y)\right).
	\end{split}\end{equation}
As in \textsc{Dalang and Nualart} \cite[(3.8)]{Dalang:04}, 
	the triangle inequality and the elementary bound 
	$(a-b)^2 \geq \frac{1}{2} a^2-b^2$ together yield
	\begin{equation}\label{biva}\begin{split}
		 p^i_{t,x;s,y}(z_1,z_2)&\leq \frac{1}{2 \pi\sigma_{s,y}\tau }
			\exp \left(- \frac{\vert z_1-z_2 \vert^2}{4 \tau^2} \right)\\
		&\hskip1.4in \times
			\exp \left(
			\frac{\vert z_2 \vert^2\, \vert 1-m \vert^2}{4 \tau^2} \right)
			\exp \left(- \frac{\vert z_2 \vert^2}{2\sigma_{s,y}^{2}} \right).
	\end{split}\end{equation}

	We will use the technical estimates in the next two lemmas in order to estimate the right-hand side of \eqref{biva}.
	
	\begin{lemma} \label{l1a}
		Fix $t_0 >0$.		There exist $c_1,c_2>0$ such that for all 
		$s,t \in [t_0\,,T]$, $x,y \in [0\,,1]$ and $i=1,\ldots,d$,
\begin{equation} \label{incr}
			\frac{1}{c_1} {\bf \Delta}((t\,,x)\,; (s\,,y)) 
				\leq \E\left[ (v_i(t\,,x)-v_i(s\,,y))^2\right] \leq c_1{\bf \Delta}
				((t\,,x)\,; (s\,,y)) 
\end{equation}			
and
\begin{equation} \label{incr2}
			\vert \sigma_{t,x}-\sigma_{s,y}\vert \leq c_2  
			  \left(\vert t - s \vert ^{1/2} + \vert x-y\vert \log\frac{1}{\vert x-y\vert} \right).
\end{equation}
	\end{lemma}

\begin{proof} 
	We assume without loss of generality that $s \leq t$.	We start by proving the upper bound in (\ref{incr}).  We note first that
$$
  \E [ (v_i(t\,,x)-v_i(s\,,y))^2 ] = \int_s^t dr \int_0^1 dz\,
  G^2_{t-r}(x\,,z)  +   \int_0^s dr \int_0^1 dz\,
  (G_{t-r}(x\,,z)-G_{s-r}(y\,,z))^2 . 
$$ 
This can be bounded above by 
		\begin{equation}\begin{split}
			& \int_s^t dr \int_0^1 dz\,  G^2_{t-r}(x\,,z) +
				2\int_0^s \int_0^1  (G_{t-r}(x\,,z)-G_{s-r}(x\,,z))^2 \, dr\, dz \\
			&\hskip2.3in +2\int_0^s \int_0^1  
				(G_{s-r}(x\,,z)-G_{s-r}(y\,,z))^2 \, dr\, dz.
		\end{split}\end{equation}
This and Lemma B.1 of \textsc{Bally, Millet, and Sanz--Sol\'e}
\cite{Bally:95} show that there is $C_0 < \infty$ such that
\begin{equation}\label{rdub1}
   \E [ (v_i(t\,,x)-v_i(s\,,y))^2 ] \leq C_0 {\bf \Delta}
	((t\,,x)\,; (s\,,y)),
\end{equation}
which is the desired upper bound.
		
		We now turn to the lower bound in (\ref{incr}). We consider three different cases.
\vskip 12pt
		
		{\em Case 1:} $s=t$. We follow 
			\textsc{Walsh} \cite[p.\ 323--326]{Walsh:86}
			and express the Green kernel for the heat equation with Neumann 
			boundary conditions as
			\begin{equation}\label{rd4.15}
				G_t(x\,,y)=\sum_{k=0}^{\infty} e^{-\pi^2 k^2 t} \phi_k(x) \phi_k(y),
			\end{equation}
			where $\phi_0(x):=1$ and $\phi_k(x):=2^{1/2} \cos(k \pi x)$ 
			[$k \geq 1$]. 
				 Therefore,
			\begin{equation}\label{rd4.16}\begin{split}
				v_i(t\,,x) &= \int_0^t \int_0^1 G_{t-s} (x\,,y)\, W^i(ds dy)\\
				&=  \sum_{k=0}^{\infty} \phi_k(x) \, A_t^k,
			\end{split}\end{equation}
			where
			\begin{equation}
				A_t^k :=\int_0^t \int_0^1 e^{-\pi^2 k^2 (t-s)} \phi_k(y) W^i(ds dy).
			\end{equation}
			We note that $t$ is fixed, and$\{A_t^k\}_{k=0}^\infty$ are 
			independent centered Gaussian random variables with variance
			\begin{equation}\begin{split}
				\text{Var}(A_t^k) &=\int_0^t ds \int_0^1 dy \, 
					e^{-2 \pi^2 k^2 s} \phi_k^2(y)\\
				&=
				\begin{cases}
					(1-e^{-2\pi^2  k^2 t})/(2\pi^2  k^2 ) & \text{if } k \geq 1, \\
					t & \text{if } k=0.
				\end{cases}
			\end{split}\end{equation}
			In fact, the $A^k$'s are Ornstein--Uhlenbeck
			processes if $k \geq 1$, and Brownian motion when $k=0$.
			Consequently, for fixed $t$,
			\begin{equation}
				v_i(t\,,x)=t^{1/2} \xi_t^0+
				\sum_{k=1}^{\infty} \phi_k(x) \, 
				\left(\frac{1-e^{-2\pi^2  k^2 t}}{2\pi^2  k^2}
				\right)^{1/2} \, \xi_t^k,
			\end{equation}
			where $\{\xi_t^k\}_{k=0}^\infty$ 
			is an i.i.d.\ sequence of standard Gaussian random variables.
			Now, recall from \textsc{Walsh}
			\cite[Exercise 3.9, p.\ 326]{Walsh:86} that
			\begin{equation}
				B_x :=x \, \xi_t^0+ \sum_{k=1}^{\infty} 
				\frac{1}{k} \, \phi_k(x) \,\xi_t^k\qquad (0 \le x \le 1)
			\end{equation}
			defines a standard Brownian motion indexed by $[0\,,1]$. Consider
			\begin{equation}\begin{split}
				R_x &:= v_i(t\,,x)-\frac{1}{\pi \sqrt{2}} B_x\\
				&=  \left(\sqrt{t}-\frac{x}{\pi \sqrt{2}}
					\right) \, \xi_t^0 + \sum_{k=1}^{\infty} 
					\phi_k(x) \, \xi_t^k\,r_k. 
			\end{split}\end{equation}
			where
			\begin{equation}
				r_k := \frac{(1-\exp(-2\pi^2 k^2 t))^{1/2}-1}{
				2^{1/2} \pi k}.
			\end{equation}
			Because $|r_k|=O(k^{-1}\exp(-2\pi^2 k^2t_0))$ as $k\to\infty$,
			$x \mapsto R_x$ is differentiable a.s., and
			\begin{equation}\begin{split}
				\E\left[ (R_x-R_y)^2 \right] & \leq 2
					\frac{\vert x-y \vert^2}{2\pi^2}+ 
					2\sum_{k=1}^{\infty} (\phi_k(x)-\phi_k(y))^2 r_k^2 \\
				&\le \vert x-y \vert^2 + 
					4 \sum_{k=1}^{\infty} \left(
					\cos(k \pi x)-\cos(k \pi y) \right)^2 r_k^2 \\
				&= \vert x-y \vert^2 + 4 
					\sum_{k=1}^{\infty} \left[ 2 
					\sin \left(k \pi \frac{x-y}{2} \right)
					\sin \left(k \pi \frac{x+y}{2} \right)
					\right]^2 r_k^2 \\
				&\leq \vert x-y \vert^2 +4 
					\sum_{k=1}^{\infty} k^2 \pi^2 \vert x-y \vert^2 r_k^2 \\
				&\leq C \vert x-y \vert^2,
			\end{split}\end{equation}
where $C$ does not depend on $t \in [t_0\,,T]$
nor on $x,y \in [0\,,1]$. It follows that
	\begin{equation}\begin{split}
				\E\left[ (v_i(t\,,x)-v_i(t\,,y))^2 \right]  &=
					\E \left[ \left(
					\frac{B_x-B_y}{\sqrt{2}} +R_x-R_y
							  \right)^2\right] \\
				&\geq \frac{1}{4} \E [(B_x-B_y)^2]- \E[(R_x-R_y)^2] \\
				&\geq \frac{1}{4} \, \vert x-y \vert - C \vert x-y \vert^2 \\
				&\geq c \, \vert x-y \vert,
			\label{rd.4.24}
	\end{split}\end{equation}
for $\vert x-y \vert$ sufficiently small and for all $t \in [t_0\,,T]$.
			
	Observe that 
\begin{equation}
	  \E\left[ (v_i(t\,,x)-v_i(t\,,y))^2\right] =
	  \int_0^t dr \int _0^1 dz \, (G_{t-r}(x\,,z) - G_{t-r}(y\,,z))^2
\end{equation}
is strictly positive,
since the integrand is not identically zero. Because
this expression is a continuous function of $(t,x,y)$,
it is bounded below on
$\{(t,x,y) \in [t_0\,,T] \times [0\,,1]^2: \vert x - y\vert \geq \varepsilon \}$
by a positive constant for every fixed $\varepsilon >0$. 
We have proved that (\ref{rd.4.24})
holds for $s=t \in [t_0\,,T]$ and $\vert x - y \vert$ sufficiently small.
Therefore, (\ref{rd.4.24}) holds for all $x,y\in [0\,,1]$ and $t \in [t_0\,,T]$
if $c$ is chosen small enough. We conclude for the moment
that there is $c >0$ such that for all $t\in[t_0\,,T]$
and $x,y \in [0\,,1]$, 
\begin{equation}\label{rdlb1}
   \E\left[ (v_i(t\,,x)-v_i(t\,,y))^2\right] \geq c \vert x - y \vert.
\end{equation}

\vskip 12pt

		{\em Case 2:} $\vert t-s\vert^{1/2} \geq \frac{c}{4C_0}
		\vert x-y\vert$, where $c$ and $C_0$ are the constants
		appearing in (\ref{rdlb1}) and (\ref{rdub1}), respectively. 
		
		By \textsc{Morien} \cite[Lemma A1.2]{Morien:99},
	\begin{equation}\label{4.26}\begin{split}
				\E\left[ (v_i(t\,,x)-v_i(s\,,y))^2 \right] 
					&\geq \int_s^t \int_0^1 G^2_{t-r}(x\,,y) \, dr\, dy\\
				&\geq \tilde c \vert t-s \vert^{1/2}.
	\end{split}\end{equation}
Because of the inequality that defines this Case 2, this is bounded below by
\begin{equation}
   \frac{\tilde c}{2}\vert t-s \vert^{1/2} + \frac{\tilde c}{2} \frac{c}{4C_0} \vert x-y\vert \geq c' {\bf \Delta} ((t\,,x)\,; (s\,,y)).
\end{equation}
This proves the lower bound in (\ref{incr}) in this Case 2.
\vskip 12pt

		{\em Case 3:} $\vert t-s\vert^{1/2} < \frac{c}{4C_0} \vert x-y\vert$,
		where $c$ and $C_0$ are the constants appearing in (\ref{rdlb1})
		and (\ref{rdub1}), respectively. 
		
		Using (\ref{rdlb1}) and (\ref{rdub1}), we observe that
\begin{equation}\begin{split}
	\E\left[ (v_i(t\,,x)-v_i(s\,,y))^2\right]
		&= \E\left[ (v_i(t\,,x) - v_i(t\,,y) + v_i(t\,,y) -v_i(s\,,y))^2\right] \\
	& \geq \frac{1}{2} \E\left[ (v_i(t\,,x)-v_i(t\,,y))^2\right]
		- \E\left[ (v_i(t\,,y)-v_i(s\,,y))^2\right]\\
	& \geq  \frac{1}{2} c \vert x-y\vert - C_0 \vert t-s\vert^{1/2}.
\end{split}\end{equation}
Because of the inequality that defines this Case 3, this is bounded below by
\begin{equation}\begin{split}
   \frac{c}{2}  \vert x-y\vert - \frac{c}{4} \vert x-y\vert
		& = \frac{c}{4} \vert x-y\vert\\
	&  \geq \frac{c}{8} \vert x-y\vert + \frac{c}{8} \frac{4C_0}{c} \vert t-s\vert^{1/2} \\
	&  \geq \min\left(\frac{c}{8}, \frac{C_0}{2}\right) {\bf \Delta} ((t\,,x)\,; (s\,,y)).
\end{split}\end{equation}
This completes the proof of Case 3 and of the lower bound in (\ref{incr}).
		
Finally we prove (\ref{incr2}). When $(t\,,x) = (s\,,y)$,
there is nothing to prove. Therefore, by the
triangle inequality, it suffices to consider the following two cases.
\vskip 12pt

		{\em (i) The case where $s=t$ and $x \neq y$.} Note that
			\begin{equation}\label{rd4.27}\begin{split}
				\vert \sigma_{t,x}-\sigma_{t,y}\vert &=
					\frac{\vert\sigma_{t,x}^2-\sigma_{t,y}^2\vert}{
					\sigma_{t,x}+\sigma_{t,y}}\\
				&\leq c \, \vert\sigma_{t,x}^2-\sigma_{t,y}^2\vert,
			\end{split}\end{equation}
		where $c$ does not depend on $t \in [t_0\,,T]$. Also, by (\ref{rd4.16}),
			\begin{equation}\begin{split}
				\sigma_{t,x}^2-\sigma_{t,y}^2 &= \sum_{k=0}^{\infty} \phi_k^2(x)
					\int_0^t ds\,  e^{-2 \pi^2 k^2 (t-s)} 
					- \sum_{k=0}^{\infty} \phi_k^2(y) \int_0^t ds  \, 
					e^{-2 \pi^2 k^2 (t-s)} \\
				&=\sum_{k=1}^{\infty} \left( \phi_k^2(x)-\phi_k^2(y)
				\right) 
					\int_0^t ds \, e^{- 2 \pi^2 k^2 s}.
			\end{split}\end{equation}
			Therefore,
			\begin{equation}\begin{split}
				\vert \sigma_{t,x}-\sigma_{t,y} \vert 
					&\leq c \, \sum_{k=1}^{\infty} 
					\frac{\vert \phi_k^2(x) - \phi_k^2(y)\vert}{k^2}\\
				&\leq 2 c \, \sum_{k=1}^{\infty}\frac{\vert 
					\phi_k(x)- \phi_k(y)\vert}{k^2}.
			\end{split}\end{equation}
			Now
			\begin{equation}\begin{split}
				\vert \phi_k(x)- \phi_k(y) \vert &\leq 4 \, 
					\left\vert \text{sin} \left(k \, \pi\,  \frac{x-y}{2}\right)
					\right\vert\\
				&\leq 4 \, \left(k \pi \, \frac{\vert x-y\vert}{2} \wedge 1 \right).
			\end{split}\end{equation}
			Consequently, as long as $|x-y|$ is sufficiently small,
			\begin{equation}\begin{split}
				\vert \sigma_{t,x}-\sigma_{t,y} \vert 
					& \leq 8c \sum_{k=1}^{\infty} \frac{1}{k^2} 
					\left(k \pi \, \frac{\vert x-y\vert}{2} \wedge 1 \right) \\
				&= \tilde{c} \left( \sum_{1\le k\le 2/\vert x-y \vert\pi}
					\frac{\vert x-y\vert}{2k}+ 
					\sum_{k> 2/\vert x-y \vert\pi}
					\frac{1}{k^2} \right) \\
				& \leq C_1 \vert x-y \vert \ln 
					\left( \frac{2}{\pi\vert x-y \vert}\right)
					+C_2|x-y|, 
			\end{split}\end{equation}
where $C_1$ and $C_2$ do not depend on $t \in [t_0\,,T]$. This proves (\ref{incr2}) when $s=t$.
\vskip 12pt

		{\em (ii) Case where $x=y$ and $s<t$.} As in (\ref{rd4.27}),
			\begin{equation}
				\vert \sigma_{t,x}-\sigma_{s,x}\vert
				\leq c \, \vert\sigma_{t,x}^2-\sigma_{s,x}^2\vert,
			\end{equation}
			and
			\begin{equation}
				\sigma_{t,x}^2-\sigma_{s,x}^2 =
				\int_s^t \int_0^1 G^2_{t-r}(x\,,y) \, dr\, dy +
				\int_0^s \int_0^1  \left( G^2_{s-r}(x\,,y)-
				G^2_{t-r}(x\,,y)\right)
				\, dr\, dy.
			\end{equation}
	We appeal to 
	\textsc{Bally, Millet, and Sanz--Sol\'e}
	\cite[Lemma B.1]{Bally:95} to see that the first term
	is bounded above in absolute value by $c(t-s)^\frac12$.
	Using (\ref{rd4.15}), we see that the second term is equal to
\begin{equation}\begin{split}
	& \sum_{k=1}^\infty \phi_k^2(x)\left(
		\int_0^s dr\,  e^{-2 \pi^2 k^2 (s-r)} 
		-  \int_0^s dr  \, 
		e^{-2 \pi^2 k^2 (t-r)}\right)\\
	&\hskip2.3in = \sum_{k=1}^\infty \phi_k^2(x)
		\left(1 - e^{-2 \pi^2 k^2 (t-s)}\right)\int_0^s dr\,
		e^{-2 \pi^2 k^2 r}.
\end{split}\end{equation}
Using the elementary inequality $0 \leq 1-e^{-x} \leq \min(x, 1)$,
valid for  all $x \geq 0$, evaluating the remaining integral
and using the fact that $\vert \phi_k^2(x) \vert \leq 2$,
we see that this is bounded above by
\begin{equation}\begin{split}
   c \sum_{k=1}^\infty \frac{\min(\pi^2 k^2 (t-s), 1)}{\pi^2 k^2} 
   &\leq C \left(\sum_{k=1}^{\pi^{-1} (t-s)^{-1/2}} (t-s)
   	+ \sum_{k > \pi^{-1} (t-s)^{-1/2}} \frac{1}{\pi^2 k^2} \right) \\
   & \leq \tilde C (t-s)^{1/2}.
\end{split}\end{equation}
This completes the proof of (\ref{incr2}) and of the lemma.
\end{proof}

	\begin{lemma} \label{l2a}
	Fix $t_0 >0$.	There exist $c_1,c_2>0$ such that for all 
		$s,t \in [t_0\,, T]$ and $x,y \in [0\,,1]$,
		\begin{align}\label{bisa}
			\frac{1}{c_1} {\bf \Delta}((t\,,x)\,; (s\,,y)) \leq 
				\sigma_{t,x}^2 \sigma^2_{s,y}- \sigma^2_{t,x;s,y }
				&\leq c_1{\bf \Delta}((t\,,x)\,; (s\,,y)),\\ \label{bisa2}
			\vert \sigma_{t,x}^2-\sigma_{t,x;s,y} \vert
				&\leq c_2\left[ {\bf \Delta}((t\,,x)\,; (s\,,y)) \right]^{1/2}.
		\end{align}
	\end{lemma}

	\begin{proof}
		Let $\gamma_{t,x;s,y}^2 := \E[(v_i(t\,,x)-v_i(s\,,y))^2]$.
		Then using \textsc{Mueller and Tribe} \cite[(4.3)]{Mueller:03},
		\begin{equation} \label{equat2}
			\sigma_{t,x}^2 \sigma^2_{s,y}- 
			\sigma^2_{t,x;s,y }= \frac{1}{4} 
			\left( \gamma_{t,x;s,y}^2-(\sigma_{t,x}-\sigma_{s,y})^2
			\right) \left( (\sigma_{t,x}+\sigma_{s,y})^2-
			\gamma_{t,x;s,y}^2 \right).
		\end{equation}
		By Lemma \ref{l1a}, 
		$\gamma_{t,x,s,y}^2 \leq c {\bf \Delta}((t\,,x)\,; (s\,,y))$.
		Therefore, the second factor of (\ref{equat2}) is bounded
		below by a positive constant when $s,t \in [t_0\,,T]$
		and $(t\,,x)$ is near $(s\,,y)$. Furthermore, another
		application of Lemma \ref{l1a} yields
		\begin{equation}\begin{split}
			\gamma_{t,x,s,y}^2-(\sigma_{t,x}-\sigma_{s,y})^2 &
				\geq c {\bf \Delta}((t\,,x)\,; (s\,,y))
				-\tilde{c}\left[ {\bf \Delta}((t\,,x)\,; (s\,,y))
				\right]^{3/2} \\
			&\geq \tilde{c} {\bf \Delta}((t\,,x)\,; (s\,,y)).
		\end{split}\end{equation}
		This proves the lower bound of (\ref{bisa}) provided $(t\,,x)$
		is sufficiently near $(s\,,y)$. 
		
In order to extend this inequality to all $(t\,,x)$ and
$(s\,,y)$ in $[t_0\,,T]\times [0\,,1]$, it suffices to show that
\begin{equation}\label{rdpos1}
		\sigma_{t,x}^2 \sigma^2_{s,y}- \sigma^2_{t,x;s,y }
		> 0 \qquad\mbox{if } (t\,,x) \neq (s\,,y).
\end{equation}
This could be proved by elementary arguments,
but since we are only interested in the conclusion,
we use results available in the literature, even if
they constitute overkill. Notice that if $s=t$ and $x \neq y$, 
then this holds because by \textsc{Bally and Pardoux}
\cite{Bally:98}, the random vector $(v_i(t\,,x),v_i(t\,,y))$
has a density with respect to Lebesgue measure. Since this
is a Gaussian random vector, this implies that the
determinant of its variance/covariance matrix is non-zero,
and this determinant is equal to
$\sigma_{t,x}^2 \sigma^2_{s,y}- \sigma^2_{t,x;s,y }$. 

   If $s < t$, and if this determinant were equal to $0$,
   then we would have $\vert \rho_{t,x;s,y}\vert = 1$,
   so there would be $\lambda \in \R$ such that
   $v_i(t\,,x) = \lambda v_i(s\,,y)$ a.s., and, in particular, we would have
\begin{equation}
 \E\left[ (v_i(t\,,x)- \lambda v_i(s\,,y))^2\right] =0.
\end{equation}
However, the left-hand side is equal to
\begin{equation}
   \int_s^t \int_0^1 G^2_{t-r}(x\,,z) \, dr\, dz +
	\int_0^s \int_0^1  \left( G_{t-r}(x\,,z)-\lambda 
	G_{s-r}(y\,,z)\right)^2 >0,
\end{equation}
which is a contradiction. Therefore, $\sigma_{t,x}^2 \sigma^2_{s,y}- \sigma^2_{t,x;s,y }>0$ when $s < t$ or $s=t$ and $x \neq y$. This completes the proof of (\ref{rdpos1}) and of the lower bound (\ref{bisa}). 

		In order to prove the upper bound of 
		(\ref{bisa}), we use Lemma \ref{l1a}, once again,
		to see that the first
		factor of (\ref{equat2}) is bounded 
		above by $c \, {\bf \Delta}((t\,,x)\,; (s\,,y))$.
		Similarly, the second factor is bounded above by a constant.
		The desired upper bound follows.

		It remains to prove (\ref{bisa2}). For this, note that
		\begin{equation}\begin{split}
			\vert \sigma_{t,x}^2-\sigma_{t,x;s,y} \vert&
				= \left\vert \gamma_{t,x;s,y}^2 + \text{Cov}
				\left( v_i(t\,,x)-v_i(s\,,y) \, , v_i(s\,,y)
				\right) \right\vert\\
			& \leq \gamma_{t,x;s,y}^2 +\gamma_{t,x;s,y} \sigma_{s,y} \\
			& \leq c \, [{\bf \Delta}((t\,,x)\,; (s\,,y))]^{1/2},
		\end{split}\end{equation}
	where we have  used Lemma \ref{l1a} twice in the last
	inequality. This implies the desired bound.
	\end{proof}

	\vskip 12pt
	By applying Lemmas \ref{l1a} and \ref{l2a} in (\ref{biva}), we 
	find, using the independence of the components
	$v_1,\dots,v_d$, that for all $z_1,z_2 \in [-M\,,M]^d$,
	\begin{equation}
		p_{t,x;s,y}(z_1\,,z_2)\leq  \frac{c}{ {\bf \Delta}((t\,,x)\,;
		(s\,,y))^{d/2}}
		\exp \left(-\frac{\Vert z_1-z_2 \Vert^2}{
		c {\bf \Delta}((t\,,x)\,; (s\,,y))} \right).
	\end{equation}
	This verifies {\bf A2}, whence follows the proof of Proposition \ref{capa}.
\end{proof}

We now establish an upper bound for hitting small balls.
Note that by Lemma \ref{l1a} and the fact that $u$ and
$v$ are Gaussian processes, Theorem \ref{eta} show that
\eqref{cond:i} holds for the solution $u$ of (\ref{e:p})
and for any $\beta \in\,]0\,,d[$. The following lemma
improves this by establishing \eqref{cond:i} for $\beta=d$,
by using the structure of the Gaussian fields $u$ and $v$.

\begin{proposition} \label{haus}
	Fix $t_0 >0$. The solution to \textnormal{(\ref{e:p})}
	satisfies \eqref{cond:i} with $\beta=d$, $I = [t_0\,,T]$ and $J = [0\,,1]$.
\end{proposition}

In order to prove Proposition \ref{haus}, we need the following lemma. 

\begin{lemma} \label{incr23}
	Let ${ u}=({ u}(t\,,x))$ be as in \textnormal{(\ref{e:p})}. 
	For all $p \geq 1$, there exists $A_p>0$ such that 
	for all $\epsilon>0$ and all $(t\,,x)$ fixed,
	\begin{equation} \label{incr4}
		\E \left[ \sup_{[{\bf \Delta}((t,x) \, ; (s,y))]^{1/2} 
		\leq \epsilon} \Vert u(t\,,x)-u(s\,,y) \Vert^p \right] 
		\leq A_p \epsilon^p.
	\end{equation}
\end{lemma}

\begin{proof}
	It suffices to prove (\ref{incr4}) for each coordinate $u_i$, $i=1,\ldots,d$.
	We plan to apply Proposition \ref{pr:Garsia} with 
	$S:=S_{\epsilon}=\{(s,y): [{\bf \Delta}((t\,,x)\,;(s\,,y))]^ {1/2} < \epsilon \}$, 
	$\rho((t\,,x)\,,(s\,,y)):=[{\bf \Delta}((t\,,x)\,;(s\,,y))]^{1/2}$, $\mu(dtdx):=dtdx$, $\Psi	(x):=e^{\vert x \vert}-1$, $p(x):=x$, and $f:=u_i$. 
	Then, by Lemma \ref{l1a} and the fact that $u = \sigma v$,
	\begin{equation}
		\E[\mathcal{C}]\le \E \left[ \int_{S_{\epsilon}} \,
		dr\, d\bar{y} \int_{S_{\epsilon}}\, ds\, dy \
		\exp\left(\frac{\vert u_i(r\,,\bar{y})-u_i(s\,,y)\vert}{
		\left(\vert r-s \vert^{1/2}+\vert \bar{y}-y \vert 
		\right)^{1/2}} \right)\right] \leq c_0 \epsilon^{12}.
	\end{equation}
	In accord with Proposition \ref{pr:Garsia}, and
	by repeated application of Jensen's inequality,
	\begin{equation}\begin{split}
		&\E \left[ \sup_{[{\bf \Delta}((t,x) \,;(s,y))]^ {1/2} \leq \epsilon}
			\vert u_i(t\,,x)-u_i(s\,,y) \vert^p \right] \\
		& \hskip1.5in \qquad \leq 8^p
			\E \left[ \left( \int_0^{2 \epsilon} du \
			\ln \left( 1+\frac{\mathcal{C}}{\left[
			\mu(B_{\rho}((t\,,x)\,,u/2) \right]^2
			} \right) \right)^p\right] \\
		& \hskip1.5in \qquad= 8^p \E 
			\left[ \left( \int_0^{2 \epsilon} du\ \ln 
			\left(1+\frac{\mathcal{C}}{c_1 u^{12}} \right) 
			\right)^p\right]   \\
		& \hskip1.5in \qquad\leq 8^p (2 \epsilon)^{p-1} 
			\E \left[  \int_0^{2 \epsilon} du\
			\ln^p \left(1+\frac{\mathcal{C}}{c_1 u^{12}} \right)
			\right] \\
		& \hskip1.5in \qquad\leq 8^p (2 \epsilon)^{p-1}
			\int_0^{2 \epsilon} du\ \ln^p \left(
			1+\frac{\E[\mathcal{C}]}{c_1 u^{12}} \right) \\
		& \hskip1.5in \qquad\leq 8^p (2 \epsilon)^{p-1}
			\int_0^{2 \epsilon} du\ \ln^p 
			\left(1+\frac{c_0}{c_1} \left(
			\frac{\epsilon}{u}\right)^{12} \right),
	\end{split}\end{equation}
	and this is manifestly a constant multiple of $\epsilon^p$.
\end{proof}


\begin{proof}[Proof of Proposition \ref{haus}]
Let ${ u}=({ u}(t\,,x))$ be as in \textnormal{(\ref{e:p})}. Let
	$R_{k,l}^n:=[t_k^n \,,t_{k+1}^n] \times [x_{\ell}^n \,,x_{\ell+1}^n]$ be as in \textnormal{(\ref{covering})}. We are going to show that there is $c < \infty$ such that for all $z \in \R^d$ and
	$\epsilon>0$,
	\begin{equation}
		\P \{ u(R_{k,l}^n) \cap B(z\,,\e) \neq \varnothing\} \le c \e^d.
	\end{equation}
	That is, $u$ satisfies \eqref{cond:i} with $\beta = d$.

Note that it suffices to prove this with $u$ replaced by $v$, where $v$ is the
solution of (\ref{v}). Without loss of generality, we set $\epsilon:=2^{-n}$. It suffices to prove that there exists 
	$c\in\,]0\,,\infty[$ such that for all $k,\ell$,
	\begin{equation}\label{StString:goal:UB}
		\P\left\{ v(R_{k,\ell}^n)\cap B(z\,,2^{-n})\neq\varnothing
		\right\} \le c \, 2^{-nd}.
	\end{equation}
	Consider
	\begin{equation}
		c^n_{k,\ell} (t\,,x) := \frac{ \E\left[ v_1(t\,,x)
		v_1(t_k^n\,,x^n_{\ell})\right]}{
		\mathrm{Var}\left[ v_1(t_k^n\,,x^n_{\ell})\right]},
	\end{equation}
	so that
	\begin{equation}
		\E\left[ v(t\,,x)\,\left|\, v(t_k^n\,,x^n_{\ell}) \right.\right] =
		c^n_{k,\ell}(t\,,x) v(t_k\,,x_{\ell}).
	\end{equation}
	Clearly,
	\begin{equation}\begin{split} \label{v2}
		\P\left\{ v(R^n_{k,\ell})\cap B(z\,,2^{-n})\neq\varnothing\right\}
			& = \P\left\{ \inf_{(t\,,x)\in R_{k,\ell}^n}
           \left\Vert v(t\,,x)-z\right\Vert\le 2^{-n} \right\}\\
		&\le \P\left\{ Y^n_{k,\ell} \le 2^{-n} + Z^n_{k,\ell}\right\},
	\end{split}\end{equation}
	where
	\begin{equation}\begin{split}
		Y_{k,\ell}^n &:= \inf_{(t,x)\in R_{k,\ell}^n} \left\Vert
			c_{k,\ell}^n(t\,,x) v(t_k^n\,,x^n_{\ell}) - z\right\Vert,\
			\text{and}\\
		Z_{k,\ell}^n &:= \sup_{(t,x)\in R^n_{k,\ell}} \left\Vert
			v(t\,,x) - c^n_{k,\ell}(t\,,x) v(t_k\,,x_{\ell})\right\Vert.
	\end{split}\end{equation}
For $r>0$,
\begin{equation}\begin{split}
		\P\{Y_{k,\ell}^n\le r\}
			&\le \P\left( \bigcap_{i=1}^d G^{i,n}_{k,\ell} \right)\\
		&= \prod_{i=1}^d \P(G^{i,n}_{k,\ell}),
\end{split}\end{equation}
where
\begin{equation}	
  G^{i,n}_{k,\ell} = \left\{\inf_{(t,x)\in R_{k,\ell}^n} \left\vert
	c_{k,\ell}^n(t\,,x) v_i(t_k^n\,,x^n_{\ell}) - z_i\right\vert \le r\right\}.
\end{equation}
The inequality $\vert c_{k,\ell}^n(t\,,x) v_i(t_k^n\,,x^n_{\ell}) - z_i\vert \le r$ is equivalent to
\begin{equation}
    \frac{z_i-r}{c^n_{k,\ell}(t,x)} \le v_i(t_k^n\,,x^n_{\ell}) \le \frac{z_i+r}{c^n_{k,\ell}(t,x)},
\end{equation}
and the interval $[(z_i-r)/c^n_{k,\ell}(t,x)\,, (z_i+r)/c^n_{k,\ell}(t,x)]$ has length bounded above by $2r/e^n_{k,\ell}$, where
\begin{equation}
		e^n_{k,\ell} := \inf_{(t,x)\in R^n_{k,\ell}}c^n_{k,\ell}(t\,,x).
\end{equation}
Therefore,
\begin{equation}
   \P(G^{i,n}_{k,\ell}) \le \sup_{x\in \R}
   \P \left\{ x\le v_i(t_k^n\,,x^n_{\ell}) \le x+\frac{2r}{e^n_{k,\ell}}
   \right\}.
\end{equation}

	Observe that for all $(t\,,x)\in R_{k,\ell}^n$,
\begin{equation}\label{eq:ac}\begin{split}
		\left| c_{k,\ell}^n(t\,,x)-1\right| &= \frac{\left|
			\E\left[ v_1(t_k^n\,,x_{\ell}^n) \cdot \left(
			v_1(t\,,x)- v_1(t_k^n\,,x_{\ell}^n)\right) \right] \right|}{
			\text{Var}\left[ v_1(t_k^n\,,x_{\ell}^n)\right]}\\
		&\le \left(\frac{
			\E\left[ \left( v_1(t\,,x)- v_1(t_k^n\,,x_{\ell}^n)\right)^2
			\right]}{\text{Var}\left[ v_1(t_k^n\,,x_{\ell}^n)
			\right]}\right)^{1/2}.
\end{split}\end{equation}
Lemma \ref{l1a} implies that the numerator is $O(2^{-n})$,
whereas the denominator is bounded below by a positive constant. Therefore,
	\begin{equation}\label{eq:c}
		\left| c_{k,\ell}^n(t\,,x) - 1\right| \le \frac{c}{2^n}\qquad
		\text{for all }(t\,,x)\in R^n_{k,\ell}.
	\end{equation}
	We emphasize the fact that the constant $c$ does not depend
	on the choice of $(n,k,\ell)$. It follows from \eqref{eq:ac}
	and \eqref{eq:c} that
\begin{equation}
	\frac{r}{e^n_{k,\ell}} \le c\, r.
\end{equation}
Since  $\{ v_i(t_k^n\,,x_{\ell}^n)\}_{i=1,\dots,d}$ are independent, centered, Gaussian random variables with variance bounded below by a positive constant,
\begin{equation}\label{eq:3}
       \P\left\{ Y_{k,\ell}^n \le r\right\} \le c\, r^d,
\end{equation}
where $c$ does not depend on our choice of $(k,\ell,n,r)$.
   Because $Y_{k,\ell}^n$ and $Z_{k,\ell}^n$ are independent,
   \eqref{v2} and \eqref{eq:3} together imply that
   \begin{equation}\label{eq:4}\begin{split}
       \P\left\{ v(R^n_{k,l}) \cap B(z\,,2^{-n})\neq\varnothing \right\}
           & \le c \, \E\left[ \left( 2^{-n}+Z^n_{k,\ell}\right)^d\right]\\
       &\le c\left( 2^{-nd} + \E\left[ (Z^n_{k,\ell})^d\right]\right).
   \end{split}\end{equation}
   We bound $Z_{k,\ell}^n$ by
   \begin{equation}
       Z_{k,\ell}^n\le Z_{k,\ell}^{(1), n} + Z_{k,\ell}^{(2), n},
   \end{equation}
   where
   \begin{equation}\begin{split}
       Z_{k,\ell}^{(1), n} & := \sup_{(t,x)\in R_{k,\ell}^n} \left\Vert
           v(t\,,x) - v(t_k^n\,,x_{\ell}^n)\right\Vert,\\
       Z_{k,\ell}^{(2),n} &:= v(t_k^n\,,x_{\ell}^n) \times
           \sup_{(t,x)\in R_{k,\ell}^n} \left|
           1 - c_{k,\ell}^n(t\,,x)\right|.
   \end{split}\end{equation}
   On one hand, \eqref{eq:c} implies that the $d$-th
   moment of $Z_{k,\ell}^{(2), n}$ is at most constant times $2^{-nd}$.
   On the other hand, Lemma \ref{incr23}
   proves that
   \begin{equation}
       \E\left[ \left( Z_{k,\ell}^{(1), n} \right)^d\right]
       \le c \, 2^{-nd}.
   \end{equation}
   Therefore, \eqref{eq:4} implies \eqref{StString:goal:UB},
   whence the proposition follows.
\end{proof}


   The main result of this section is the following theorem, which summarizes the preceding results.

\begin{theorem}\label{rd:thm4.9}
	Let $u=(u(t\,,x))_{t \in [0,T], x\in [0\,,1]}$ be the
	solution of (\ref{e:p}). Fix $t_0 >0$. Then the
	conclusions of Theorems \ref{thm:cap:LB}, \ref{thm:cap2:LB},
	\ref{thm:meas:UB}, and \ref{thm:cap2:LBB} hold for $u$,
	with $I=[t_0\,,T]$, $J=[0\,,1]$, and $\beta = d$.
\end{theorem}

\begin{proof}
By Proposition \ref{capa}, A.1 and A.2 are satisfied
for $u$ with these choices of $I$, $J$ and $\beta$.
Therefore, the conclusions of Theorems \ref{thm:cap:LB},
\ref{thm:cap2:LB} are also satisfied. By Proposition \ref{haus},
$u$ satisfies (\ref{cond:i}) with $\beta = d$, $I=[t_0\,,T]$ and
$J=[0\,,1]$. Therefore, the conclusions of Theorems
\ref{thm:meas:UB}, and \ref{thm:cap2:LBB} are also satisfied.
\end{proof}

\begin{remark} We could have considered the system
(\ref{eq:spde}) with Dirichlet boundary conditions
instead of the Neumann boundary conditions (\ref{neumann}).
In this case, the results and proofs are essentially unchanged,
except that one must replace the interval $J=[0\,,1]$
by $J=[\e\,,1-\e]$, where $\e >0$ is fixed. Indeed,
a lower bound such as (\ref{4.26}) would obviously
not be satisfied at $x=y=0$ or $x=y=1$ with Dirichlet boundary conditions.
\end{remark}

\section{The case of additive noise}\label{sec5}

The aim of this section is to transfer the results 
of Section 4 for the Gaussian process (\ref{e:p})
to the process (\ref{spde}). Subsequently,
we will establish Theorem \ref{th:hitprob}
and Corollary \ref{cor:hitprob} of the Introduction.
For this, we will use the 
following general fact which is a consequence of Girsanov's
theorem.

\begin{proposition} \label{fact}
	Let ${ u}$ denote the solution of \textnormal{(\ref{eq:spde})} 
	and let ${v}$ denote the solution of \textnormal{(\ref{eq:spde})} 
	with ${b}\equiv 0$, that is, ${v}$ is the the solution of 
	\textnormal{(\ref{e:p})}. Then for any $\epsilon>0$, there exists $c>0$ such that for all be a Borel subsets $B$ of $C([0\,,T] \times [0\,,1],\, \R^d)$, 
	\begin{equation} \label{bounds}
		\frac 1c (\P\{{ v} \in B\})^{1+\epsilon} 
		\leq \P\{{ u} \in B\} \leq c \, 
		(\P\{{ v} \in B\})^{1/(1+\epsilon)}.
	\end{equation}
\end{proposition}

\begin{proof}
	We follow the proof of Corollary 5.3 of \textsc{Dalang and Nualart}
	\cite{Dalang:04} and consider
	\begin{equation}\begin{split}
	L_t &:= \exp \biggl(-\int_0^t \int_0^1 G_{t-s}(x\,,y) 
		\,{{\sigma}}^{-1} { b}(u(s\,,y)) \cdot W(ds \, dy) \\
	&\hskip2in -
		\frac{1}{2} \int_0^t \int_0^1 (G_{t-s}(x\,,y))^2 \, 
		\Vert {{\sigma}}^{-1} { b}(u(s\,,y))\Vert^2 \, ds\, dy\biggr), \\
	J_t &:= \exp \biggl(-\int_0^t \int_0^1 G_{t-s}(x\,,y) 
		\, {{\sigma}}^{-1} { b}(v(s\,,y)) \cdot W(ds\, dy) \\
	&\hskip2in +
		\frac{1}{2} \int_0^t \int_0^1 (G_{t-s}(x\,,y))^2 
		\, \Vert {{\sigma}}^{-1} { b}(v(s\,,y))\Vert^2 \, ds
		\, dy\biggr).
	\end{split}\end{equation}
	Let ${\rm Q}$ denote the probability measure defined by
	\begin{equation}
		\frac{d {\rm Q}}{d\P}(\omega)=L_t(\omega).
	\end{equation}
	Then, by Girsanov's theorem as stated in Proposition
	1.6 of \textsc{Nualart and Pardoux} \cite{NP} 
	(see also \textsc{Dalang and Nualart} \cite[Theorem 5.2]{Dalang:04}),
	\begin{equation} \label{girsanov}
		\P\{{ u} \in B\} = \E_\P
		\left[ \1_{\{{ u} \in B\}}\right]=
		\E_{{\rm Q}} \left[ \1_{\{{ u} \in B\}} L_t^{-1}
		\right]=\E_\P\left[\1_{\{{ v} \in B\}} J_t^{-1}\right].
	\end{equation}
	Let $\epsilon>0$ and apply H\"older's inequality to find that
\begin{equation}\begin{split}
   \P\{v\in B\} &= \E_\P \left[ \1_{\{{ v} \in B\}}
   	J_t^{-1/(1+\e)}\, J_t^{1/(1+\e)} \right]\\ 
   &\le \left(\E_\P\left[\1_{\{{ v} \in B\}} J_t^{-1}\right]
   	\right)^{1/(1+\e)} \left(\E_\P\left[ J_t^{1/\e}\right]\right)^{\e/(1+\e)},
\end{split}\end{equation}
and therefore,
	\begin{equation}
		\P\{{ u} \in B\} \ge (\P\{{ v} \in B\} )^{1+\epsilon} 
		\left( \E_\P\left[ J_t^{1/\epsilon}
		\right] \right)^{-\epsilon}.
	\end{equation}
	Let $r=1/\e$. By the Cauchy--Schwarz inequality,
	\begin{equation}\begin{split}
		\E_\P[J_t^r] &\leq \biggl( \E_\P \biggl[
			\exp \biggl( \int_0^t \int_0^1 -2 r \, G_{t-s}(x\,,y)
			\,{{\sigma}}^{-1} { b}(v(s\,,y)) \cdot W(ds \, dy) \\
		&\hskip1in  -\frac12
			\int_0^t \int_0^1 4 r^2 (G_{t-s}(x\,,y))^2 \,
			\Vert {{\sigma}}^{-1} { b}(v(s\,,y))\Vert^2 
			\, ds\, dy \biggr) \biggr] \biggr)^{1/2} \\
		& \quad \times \biggl(
			\E_\P \biggl[ \exp \biggl( \int_0^t \int_0^1 (2r^2+r) 
			(G_{t-s}(x\,,y))^2 \, \Vert {{\sigma}}^{-1} 
			{ b}(v(s\,,y))\Vert^2 \, ds\, dy \biggr) 
			\biggr] \biggr)^{1/2}.
	\end{split}\end{equation}
	The first expectation on the right-hand side 
	equals $1$ since it is the expectation of 
	an exponential martingale with
	bounded quadratic variation. The second factor 
	is bounded by some positive finite constant.
	This proves the lower bound of (\ref{bounds}).

	In order to prove the upper bound, let $\epsilon>0$ and
	apply H\"older's inequality to the right-hand side of (\ref{girsanov}):
	\begin{equation}
		\P\{{ u} \in B\} \le (\P\{{ v} \in B\} )^{1/1+\epsilon} 
		\left( \E_\P\left[ J_t^{-(1+\epsilon)/\epsilon} \right]
		\right)^{\epsilon/1+\epsilon}.
	\end{equation}
	Let $r= (1+\e)/\e$. Again by the Cauchy-Schwarz inequality,
	\begin{equation}\begin{split}
		\E_\P[J_t^{-r}] &\leq \biggl( \E_\P 
			\biggl[\exp \biggl(\int_0^t 
			\int_0^1 2 r \, G_{t-s}(x\,,y) 
			\,{{\sigma}}^{-1} { b}(v(s\,,y)) \cdot W(ds \, dy) \\
		&\hskip1in  -\frac12
			\int_0^t \int_0^1 4 r^2 (G_{t-s}(x\,,y))^2
			\, \Vert {{\sigma}}^{-1} { b}(v(s\,,y))\Vert^2
			\, ds\, dy \biggr) \biggr] \biggr)^{1/2} \\
		& \quad\times \biggl(
			\E_\P \biggl[\exp \biggl( \int_0^t \int_0^1 (2r^2-r) 
			(G_{t-s}(x\,,y))^2 \, \Vert {{\sigma}}^{-1} 
			{ b}(v(s\,,y))\Vert^2 \, ds\, dy \biggr) \biggr]
			\biggr)^{1/2}.
	\end{split}\end{equation}
	As above, the first expectation on the right-hand 
	side equals $1$ since it is the expectation of an 
	exponential martingale with
	bounded quadratic variation and the second factor is 
	bounded above by some positive finite constant.
	This concludes the proof.
\end{proof}

Theorem \ref{th:hitprob} will be a consequence of our
next result.

\begin{proposition} \label{final}
	Let ${u}$ denote the solution of \textnormal{(\ref{eq:spde})}. 
	Let $I \subset\,]0\,,T]$ and $J \subset [0\,,1]$ 
	be two fixed non-trivial compact intervals. Fix $M>0$.
	\begin{enumerate}
		\item[\textnormal{(1)}] For any $\epsilon>0$, there exists $c>0$ 
			such that for all Borel sets $A\subseteq[-M\,,M]^d$,
			\[
				\frac 1c (\Cap_{d-6}(A))^{1+\epsilon} \le
				\P\left\{ { u}(I\times J) \cap A \neq\varnothing \right\}  \le c \, (\mathcal{H}_{
	            d-6}(A))^{1/(1+\epsilon)}.
	   \]
		\item[\textnormal{(2)}] For all $t \in\,]0\,,T]$	and $\epsilon>0$, there exists $c>0$	such that for all Borel sets $A\subseteq[-M\,,M]^d$,
	        \[
				\frac 1c (\Cap_{d-2}(A))^{1+\epsilon}  
				\le \P\left\{ { u}(\{t\}\times J) \cap A
				\neq\varnothing \right\} 
				\le c \,  (\mathcal{H}_{
	            d-2}(A))^{1/(1+\epsilon)}.
	\]
		\item[\textnormal{(3)}] For all $x \in [0\,,1]$ 
			and $\epsilon>0$, there exists $c>0$ such that for all Borel sets
	        $A\subseteq[-M\,,M]^d$,
	        \[
		        \frac1c (\Cap_{d-4}(A))^{1+\epsilon}
		        \le \P\left\{ { u}(I\times \{x\}) \cap A
		        \neq\varnothing \right\} \le c \, 
		        (\mathcal{H}_{d-4}(A))^{1/(1+\epsilon)}.
	 \]
	    \end{enumerate}
\end{proposition}

\begin{proof}
	In order to prove the upper bound in (1), we apply
	Proposition \ref{fact} with
	$B=\{f\in C([0\,,T]\times [0\,,1],\R^d): f(I\times J)
	\cap A \neq\varnothing \}$ and then use Theorem \ref{rd:thm4.9}. 
	When $A$ is compact, we get the lower bound in (1) in the same way. 
	Now consider the case where $A$ is a Borel set. 
	We recall that $\Cap_{\beta}$ is a Choquet capacity;
	see \textsc{Dellacherie and Meyer} \cite[Chapter 3]{DM}.
	In particular, for any Borel set $A$,
	\begin{equation}\label{choq}
		\sup_{F \subset A, \,  F \text{ compact}} \Cap_\beta(F)=\Cap_\beta(A).
	\end{equation}
	Therefore, if $F \subset A$ is compact, then
	\begin{equation}
	   \P\left\{ { u}(I\times J) \cap A \neq\varnothing \right\} \ge
	   \P\left\{ { u}(I\times J) \cap F \neq\varnothing \right\} \ge
	   \frac 1c (\Cap_{d-6}(F))^{1+\epsilon}.
	\end{equation}
	Taking, on the right-hand side, the supremum
	over such $F$ and using (\ref{choq}) proves the
	lower bound in (1) for $A$.

	The proofs of (2) and (3) are similar and are left to the reader.
\end{proof}

We now prove Theorem \ref{th:hitprob}.

\begin{proof}[Proof of Theorem \ref{th:hitprob}]
This theorem is an immediate consequence of Proposition \ref{final}.
\end{proof}


We prove Corollary \ref{cor:hitprob} next.

\begin{proof}[Proof of Corollary \ref{cor:hitprob}]
	We first prove (a). Let $z \in \mathbb{R}^d$. 
	If $d <6$, then
	$\textnormal{Cap}_{d-6}(\{ z\})=1$. Hence, the 
	lower bound of Proposition \ref{final}(1) implies that
	$\{ z\}$ is not polar. On the other hand, if $d>6$, then
	$\mathcal{H}_{d-6}(\{ z\})=0$ and the upper bound of 
	Proposition \ref{final}(1) implies that
	$\{ z\}$ is polar. 
	If $d=6$, we observe that {\sc Mueller and Tribe}
	\cite[Corollary 4]{Mueller:03} show that the law
	of their {\em stationary pinned 
	string} \cite[(2.1)]{Mueller:03} is mutually equivalent,
	on compact subsets of $]0\,,T[ \times ]0\,,1[$,
	to the law of the solution of (\ref{eq:spde})
	(see \cite[Corollary 4]{Mueller:03}).
	In this corollary, Mueller and Tribe consider
	the heat equation on the circle instead of the
	heat equation on $[0\,,1]$; however, 
	the Green's functions of these two equations are not
	very different and the proofs of \cite{Mueller:03}
	apply essentially without changes to our setting. Therefore,
	from \cite[Theorem 1]{Mueller:03} and Proposition \ref{final},
	we conclude that when $d=6$, a.s., the solution of
	(\ref{eq:spde}) does not hit points. This proves (a). 
	
	For (b), the cases $d<2$ and $d>2$ are proved exactly
	along the same lines using Proposition \ref{final}(2).
	For the case $d=2$, we again use the mutual equivalence
	of our process with the stationary pinned string of
	\cite{Mueller:03}. For $t$ fixed, the stationary
	pinned string as a function of $x$ has the same
	increments as those of a standard Brownian motion with
	values in $\R^d$ \cite[Section 2]{Mueller:03}.
	Therefore, points are polar for $x \mapsto u(t\,,x)$
	when $d=2$. This proves (b).
	
	 For (c), the statement only concerns the cases
	 $d<4$ and $d>4$, which are proved as above using Proposition \ref{final}(3).
\end{proof}

The following is another consequence of our work.

\begin{corollary} \label{codim2}
	Let ${u}$ denote the solution of \textnormal{(\ref{eq:spde})}.
	\begin{itemize}
		\item[\textnormal{(a)}] If $d\ge 6$, then
			$\dimh  
			({ u}(]0\,,T] \times\,]0\,,1[)) = 6$
			a.s.
		\item[\textnormal{(b)}] Fix $t \in\,]0\,,T]$. 
			If $d\ge 2$, then
			$\dimh  ({ u}(\{t\} \times\,]0\,,1[)) = 2$
			a.s.
		\item[\textnormal{(c)}] Fix $x \in\,]0\,,1[$. If $d\ge 4$, then
			$\dimh  ({ u}(\mathbb{R}_+ \times \{x\})) = 4$
			a.s.
	\end{itemize}
\end{corollary}
In the special case that $b_i\equiv 0$
and $\sigma_{i,j}\equiv \delta_{i,j}$, 
\textsc{Wu and Xiao} \cite{WX:07}
find a connection between \eqref{spde}
and the theory of local non-determinism,
and hence deduce Corollary \ref{codim2};
see their Theorem 2.3 and Proposition 2.4 (\emph{loc.\ cit.}).
Presently, we use an indirect and elementary codimension
argument to achieve a similar effect for the more general functions
$b_i$ and $\sigma_{i,j}$ under consideration here.

\begin{proof}
	Let $E$ be a random set. When it exists,
	the codimension of $E$ is the real number 
	$\beta \in [0\,,d]$ such that for all compact sets $A \subset
	\mathbb{R}^d$,
	\begin{equation}
		\P \{ E \cap A \neq \varnothing \}
		\begin{cases}
			>0 & \text{whenever $\dimh(A)> \beta$}, \\
			=0 & \text{whenever $\dimh(A)< \beta$}.
		\end{cases}
	\end{equation}
	See \textsc{Khoshnevisan}
	\cite[Chap.11, Section 4]{Khoshnevisan:02}. When it is
	well defined, we write
	the said codimension as $\text{codim}(E)$.
	Proposition \ref{final} implies that for $d \geq 1$:
	$\textnormal{codim}(u(\mathbb{R}_+ \times\,]0\,,1[))=(d-6)^+$;
	$\textnormal{codim}(u(\{t\} \times\,]0\,,1[))=(d-2)^+$; and
	$\textnormal{codim}(u(\mathbb{R}_+ \times \{x\}))=(d-4)^+$. 
	According to Theorem 4.7.1 of
	\textsc{Khoshnevisan}
	\textnormal{\cite[Chapter 11]{Khoshnevisan:02}},
	given a random set $E$ in $\mathbb{R}^d$ whose
	codimension is strictly between $0$ and $d$,
	\begin{equation} \label{codim}
		\dimh E+\textnormal{codim } E=d
		\qquad\text{a.s.\ on}\ \{E\neq\varnothing\}.
	\end{equation}
	 When $d>6$, this implies (a). When $d>2$, this implies (b),
	 and when $d>4$ this implies (c) of the corollary.
	 
	 For the remaining ``critical cases'' we consider the case
	 $d=6$ and prove (a) only. The corresponding results 
	 for (b) ($d=2$) and (c) ($d=4$) are proved analogously.
	 
	 Because $d=6$, it follows immediately that
	 the Hausdorff dimension of $u(]0\,,T]\times\,]0\,,1[)$ is
	 at most $6$. For the lower bound, we note that
	 $u(]0\,,T]\times\,]0\,,1[)$ will hit $A\subset\R^6$ as long as
	 $A$ has positive logarithmic capacity (Proposition \ref{final}).
	 In particular, the codimension of $u(]0\,,T]\times\,]0\,,1[)$ is
	 zero. 
	 
	 Choose and fix $\beta\in\,]0\,,6[$.
	 By Peres's Lemma (\textsc{Khoshnevisan} \cite[p.\ 436]{Khoshnevisan:02}),
	 we can find an independent
	 closed random set $\Lambda_\beta\subset\R^6$
	 such that for all $\sigma$-compact sets
	 $E\subset\R^6$: (i) $\dimh \Lambda_\beta\cap E=
	 \dimh E - \beta$ a.s.; 
	 (ii) $\P\{ \Lambda_\beta\cap E=\varnothing\}=1$
	 if $\dimh E<\beta$; and (iii) $\P\{\Lambda_\beta \cap E\neq
	 \varnothing\}\in\{0\,,1\}$.
	 Because $\dimh \Lambda_\beta=6-\beta$ is positive,
	 $\Lambda_\beta$ has positive logarithmic capacity;
	 this follows from Frostman's theorem \textsc{Khoshnevisan} \cite[p.\ 521]{Khoshnevisan:02}.
	 Therefore, by Proposition \ref{final} and (iii),
	 $u(]0\,,T]\times\,]0\,,1[) \cap \Lambda_\beta \neq\varnothing$
	 a.s. But thanks to (ii),
	 $\dimh u(]0\,,T]\times\,]0\,,1[)
	 \ge \beta$. Let $\beta\uparrow 6$
	 to deduce (a) in the case that $d=6$. This concludes
	 the proof.
\end{proof}

\begin{proposition} \label{final2}
Let ${u}$ denote the solution of \textnormal{(\ref{eq:spde})}. 
Then for all $\epsilon>0$ and $R>0$, 
there exists a positive and finite constant $a$ such that
the following holds for all compact sets 
$E\subset\,]0\,,T]\times\,]0\,,1[$,
$F\subset\,]0\,,T]$, and $G\subset\,]0\,,1[$, 
and for all $z\in B(0\,,R)$:
    \begin{enumerate}
        \item[\textnormal{(1)}]
			$a^{-1} (\Cap_{d/2}^{{\bf \Delta}}
			(E))^{1+\e} \le \P\{ \mathcal{L}(z\,;u) \cap E
			\neq\varnothing \} \le  a \, 
			(\mathcal{H}^{{\bf \Delta}}_{d/2}(E))^{1/(1+\epsilon)}$;
        \item[\textnormal{(2)}]
			$a^{-1} (\Cap_{(d-2)/4}(F))^{1+\e} 
			\le \P\{ \mathcal{T}(z\,;u)\cap 
			F\neq\varnothing\}\le  a \, 
			(\mathcal{H}_{(d-2)/4}(F))^{1/(1+\epsilon)}$;
        \item[\textnormal{(3)}]
			$a^{-1}  (\Cap_{(d-4)/2} (G))^{1+\e} 
			\le\P\{\mathcal{X}(z\,;u)\cap G\neq\varnothing\}
			\le  a \, (\mathcal{H}_{(d-4)/2}(G))^{1/(1+\epsilon)}$;
        \item[\textnormal{(4)}] For all $x \in\,]0\,,1[$,
            $a^{-1} (\Cap_{d/4}(F))^{1+\e} \le 
            \P\{\mathcal{L}_x(z\,;u) \cap F\neq
            \varnothing\}\le  a \, (\mathcal{H}_{d/4}(F))^{1/(1+\epsilon)}$;
        \item[\textnormal{(5)}] For all $t\in\,]0\,,T]$,
			$a^{-1} (\Cap_{d/2}(G))^{1+\e} \le 
			\P\{ \mathcal{L}^t(z\,;u) \cap G\neq\varnothing\}
			\le  a \, (\mathcal{H}_{d/2}(G))^{1/(1+\epsilon)}$.
    \end{enumerate}
\end{proposition}

\begin{proof}
	In order to prove (1), it suffices to use 
	Proposition \ref{fact} with 
	$B=\{f: \mathcal{L}(z\,;u) \cap E \neq\varnothing \}$
	and apply Theorem \ref{rd:thm4.9}. 
	The proofs of (2)--(5) follow in exactly the same way.
\end{proof}

\begin{corollary} \label{c5}
Let ${u}$ denote the solution of \textnormal{(\ref{eq:spde})}. 
	Choose and fix $z\in\R^d$.
	\begin{itemize}
		\item[\textnormal{(a)}] If $2\le d< 6$, then
			$\dimh\, \cT(z \, ; u) = \frac14 (6-d)$
			a.s.\ on $\{\cT(z \, ; u)\neq\varnothing\}$.
		\item[\textnormal{(b)}] If $4\le d< 6$, then
			$\dimh\, \cX(z \, ; u) = \frac12 (6-d)$
			a.s.\ on $\{\cX(z \, ; u)\neq\varnothing\}$.
		\item[\textnormal{(c)}] If $1\leq d< 4$, then
			$\dimh\, \cL_x(z \, ; u) = \frac14 (4-d)$
			a.s.\ on $\{\cL_x(z \, ; u)\neq\varnothing\}$.
		\item[\textnormal{(d)}] If $d=1$, then
			$\dimh\, \cL^t(z\, ; u) = \frac12(2-d)=\frac{1}{2}$
			a.s.\ on $\{\cL^t(z \, ;u)\neq\varnothing\}$.
	\end{itemize}
	In addition, all four right-most events have positive probability.
\end{corollary}

\begin{proof}
	The final positive-probability assertion is an
	  immediate consequence of Proposition \ref{final2} and
	Taylor's theorem \textsc{Khoshnevisan} \cite[Corollary 2.3.1 p.~523]{Khoshnevisan:02}. 
	
	For the remainder of the corollary,
	we proceed as we did in the proof of Corollary \ref{codim2}.
	By Proposition \ref{final2}, for $d \geq 1$, it holds that
	$\textnormal{codim}(\cT(z \, ; u))=\frac14 (d-2)^+$,
	$\textnormal{codim}(\cX(z \, ; u))=\frac12 (d-4)^+$,
	$\textnormal{codim}(\cL_x(z ; \, u))=d/4$, $\textnormal{codim}(\cL^t(z \, ; u))=d/2$. 
	Hence, (\ref{codim}) gives the desired statements of the corollary
	 in all but the critical cases. The critical cases are handled 
	 as was done in the proof of Corollary \ref{codim2}.
\end{proof}

\begin{remark} It is natural to expect that if $1 \leq d < 6$, then the  $\mathcal{H}^{{\bf \Delta}}$-Hausdorff dimension of $\mathcal{L}(z\,;u)$ is $(6-d)/2$. Indeed, since the $\mathcal{H}^{{\bf \Delta}}$-Hausdorff dimension of $]0,T] \times [0,1]$ is $3$, this would be compatible with the codimension argument, if it applied.
\end{remark}

\appendix

\section{Appendix: An anisotropic Kolmogorov Continuity Theorem}\label{apendB}

We first present an improvement of the classical lemma of
\textsc{Garsia} \cite{Garsia:72}.
Recall that $\Psi:\R\to\R_+$ is a \emph{strong Young function} if it is even and convex on $\R$, and strictly increasing on $\R_+$. Its {\em inverse} is $\Psi^{-1}: \R_+ \to \R_+$.

\begin{proposition}\label{pr:Garsia}
    Let $(S\,,\rho)$ be a metric space, $\mu$ a Radon measure
    on $S$,
    and $\Psi:\R\to\R_+$ a strong Young function
    with $\Psi(0)=0$ and $\Psi(\infty)=\infty$. Suppose $p:[0\,,\infty[\,
    \to \R_+$ is continuous, even, and strictly
    increasing on $\R_+$, with
    $p(0)=0$. Define, for any continuous function
    $f: S\to\R$,
    \begin{equation}
        \mathcal{C} := \iint \Psi\left(
        \frac{f(x)-f(y)}{p(\rho(x\,,y))}\right)\, \mu(dx)\, \mu(dy).
    \end{equation}
    Let $B_{\rho}(s\,,r)$ denote the open
    $d$-ball of radius $r>0$ about $s\in S$. 
    Then, for all $s,t\in S$,
    \begin{equation}\label{eq:garsia2}\begin{split}
        &|f(t)-f(s)|\\
        &\quad \le 4 \int_0^{2 \rho(s,t)}\left[
            \Psi^{-1}\left(\frac{\mathcal{C}}{\left[ \mu(B_{\rho}(s\,,u/2)\right]^2}
            \right) + \Psi^{-1}\left(\frac{\mathcal{C}}{
            \left[ \mu(B_{\rho}(t\,,u/2)\right]^2} \right)
            \right] \, p(du).
    \end{split}\end{equation}
\end{proposition}

\begin{remark}
	\textnormal{(a)}
		The following ``majorizing-measure condition'' is 
		a ready but useful
		consequence: If $\mathcal{C} < \infty$ then for all $\e>0$,
		\begin{equation}
			\sup_{s,t \in S:\ \rho(s,t) \leq \e}
			\vert f(t)-f(s) \vert \le 8 \sup_{x \in S} \int_0^{2 \e}
			\Psi^{-1}\left(\frac{\mathcal{C}}{\left[ \mu(B_{\rho}(x\,,u/2)\right]^2}
			\right)\, p(du).
		\end{equation} 
		Another extension is found in \textsc{Arnold and Imkeller}
		\cite{Arnold-Imkeller}.
		
	\textnormal{(b)}
		Suppose, instead of continuity, 
		 that $f\in L^1_{\text{\it loc}}(\mu)$ and
		\begin{equation}\label{cond:leb}
			\lim_{\e\to 0^+} \frac{1}{\mu(B_\rho(x\,,\e))}
			\int_{B_\rho(x\,,\e)} f\,d\mu = f(x)\qquad
			\text{for $\mu$-almost all $x$}.
		\end{equation}
		Then, a straight-forward modification of our proof
		shows that there is a $\mu$-null set $N$ such that
		\eqref{eq:garsia2} holds for all $s\,,t \in S\setminus N$.
		
	\textnormal{(c)} This proposition implies 
	various known Poincar\'e inequalities and 
	Besov--Morrey--Sobolev embedding theorems in metric spaces.
	A portion of this assertion in proved in 
	\textsc{Kassmann} \cite{Kassmann}
	who uses the inequality of \textsc{Arnold and Imkeller}
	\cite{Arnold-Imkeller}
	instead of ours. \textsc{Buckley and Koskela}
	\cite{BuckleyKoskelaII,BuckleyKoskelaI}
	contain some of the recent work on Sobolev embedding theory.
\label{remgarsia}
\end{remark}

\begin{proof}
    Throughout, we choose and fix $s,t\in S$, and
    follow the ideas of \textsc{Garsia} \cite{Garsia:72} closely.
    We may, and will, assume without loss of generality
    that $\mathcal{C}<\infty$. Otherwise, there is nothing to prove
    because $\Psi(\infty)=\infty$.

    Define, for any bounded set $Q\subset S$ with $\mu(Q)>0$,
    \begin{equation}
        \bar{f}_Q := \frac{1}{\mu(Q)} \int_Q f\,d\mu.
    \end{equation}
    Let $r_0:=\rho(s\,,t)$, and $Q_{-1}$ be the open $\rho$-ball centered at $s$
     of radius $r_{-1}:=2r_0$, so that
    $B_{\rho}(s\,,r_0)\cup B_{\rho}(t\,,r_0) \subset Q_{-1}$. Then, define $r_n$ iteratively by
    $p(2r_n)=\frac12 p(2r_{n-1})$ for all $n\ge 1$.
    Notice that as $n$ tends to infinity, both $r_n$ and $p(2r_n)$ decrease (by induction) to zero (by contradiction: if $\inf r_n >0$, then $p(2 \inf r_n) = \frac12 p(2\inf r_n)$, and therefore $\inf r_n =0$). 

    Define $Q_n:= B_{\rho}(s\,,r_n)$ for all
    $n\ge 0$, and apply Jensen's
    inequality to find that
 \begin{align}\nonumber
        \Psi\left( \frac{\bar{f}_{Q_n} -
		\bar{f}_{Q_{n-1}} }{p(2r_{n-1})} \right) 
		&= \Psi \left(\frac{1}{\mu(Q_n)\cdot \mu(Q_{n-1})} \int_{Q_n}\mu(dx) \int_{Q_{n-1}}\mu(dy)\, \frac{f(x)-f(y)}{p(2r_{n-1})} \right) \\
		& \le \frac{1}{\mu(Q_n)\cdot \mu(Q_{n-1})}
		\int_{Q_n}\mu(dx) \int_{Q_{n-1}}\mu(dy)\,  \Psi\left(
		\frac{f(x)-f(y)}{p(2r_{n-1})}\right).
\end{align}
    If $x\in Q_n$ and $y\in Q_{n-1}$, then
    $\rho(x\,,y) \le 2r_{n-1}$, whence $p(\rho(x\,,y))\le p(2r_{n-1})$.
    Therefore,
    \begin{equation}\begin{split}
        \Psi\left( \frac{\bar{f}_{Q_n} -
        \bar{f}_{Q_{n-1}} }{p(2r_{n-1})} \right)
        \le \frac{\mathcal{C}}{\left[ \mu(Q_n)\right]^2}.
    \end{split}\end{equation}
    Equivalently,
    \begin{equation}
        \left| \bar{f}_{Q_n} - \bar{f}_{Q_{n-1}} \right|
        \le \Psi^{-1}\left( \frac{\mathcal{C}}{\left[ \mu(Q_n) \right]^2}\right)
        p(2r_{n-1}).
    \end{equation}
    Because $p(2r_n)-p(2r_{n+1})= \frac14 p(2r_{n-1})$,
    \begin{equation}
        \left| \bar{f}_{Q_n} - \bar{f}_{Q_{n-1}} \right|
        \le 4\Psi^{-1}\left( \frac{\mathcal{C}}{\left[ \mu(Q_n) \right]^2}\right)
        \left[ p(2r_n)-p(2r_{n+1})\right].
    \end{equation}
    Note that $\cap_{n=1}^\infty Q_n=\{s\}$, whence
    $\lim_{n\to\infty}\bar{f}_{Q_n}= f(s)$ by continuity. Therefore,
    we can add the preceding over all $n\ge 1$ to find that
    \begin{equation}\begin{split}
        \left| f(s) - \bar{f}_{Q_{-1}} \right|
            &\le 4\sum_{n=0}^\infty
            \Psi^{-1}\left( \frac{\mathcal{C}}{\left[ \mu(Q_n) \right]^2}\right)
            \left[ p(2r_n)-p(2r_{n+1})\right]\\
        &\le 4\int_0^{2r_0} \Psi^{-1}\left(
            \frac{\mathcal{C}}{\left[ \mu(B_{\rho}(s\,,u/2))\right]^2
            }\right)\, p(du).
    \end{split}\end{equation}
    The same bound holds if we replace $s$ by $t$ throughout.
    This is because $t\in Q_{-1}$ as well. Therefore,
    $|f(s)-f(t)|$ is bounded above by
    \begin{equation}
		4\int_0^{2r_1}\Psi^{-1}\left(
		\frac{\mathcal{C}}{
		\left[ \mu(B_{\rho}(s\,,u/2))\right]^2 }\right)\, p(du)
		+ 4 \int_0^{2r_1}
		\Psi^{-1}\left(\frac{\mathcal{C}}{
		\left[ \mu(B_{\rho}(t\,,u/2)\right]^2}
		\right)\, p(du).
    \end{equation}
    We obtain the proposition by recalling merely that
    $r_1\le r_0 := \rho(s\,,t)$.
\end{proof}



\begin{corollary}\label{AKCT}
	Choose and fix two nonrandom compact intervals $I\subset \R$ and
	$J\subset \R$, and let $\{ v(t\,,x)\}_{t\in I, x\in J}$ 
	denote a real-valued stochastic process. Suppose that
	there exist finite constants $p> 1$, $q>0$, and $c>0$
	such that for all $(t\,,x) \in I\times J$ and $(s\,,y) \in I\times J$,
	\begin{equation}\label{rd:momentb}
	   \E(\vert v(t\,,x) - v(s\,,y)\vert^p) \leq
	   c [{\bf \Delta}((t\,,x)\,;(s\,,y))]^{3+q}.
	\end{equation}
	Then $v$ has a continuous version $\tilde v$, and for any $\alpha \in [0\,,q/p[$,
	\begin{equation}\label{A.12}
	   \E\left[\left(\sup_{(t,x) \neq (s,y)}
	   \frac{\vert \tilde v(t\,,x) - \tilde v(s\,,y)\vert}{
	   [{\bf \Delta}((t\,,x)\,;(s\,,y))]^{\alpha}} \right)^p \right] < \infty.
	\end{equation}
	In particular, there is a non-negative random variable $C$
	with $\E[C] < \infty$ such that a.s.,
	\begin{equation}
	   \vert \tilde v(t\,,x) - \tilde v(s\,,y)\vert
	   \le C [{\bf \Delta}((t\,,x)\,;(s\,,y))]^{\alpha}.
	\end{equation}

\end{corollary}

\begin{proof}
We observe that (\ref{rd:momentb}) imples that $(t\,,x) \mapsto v(t\,,x)$ 
is continuous in probability, and therefore, has a
measurable version 
(\textsc{Dellacherie and Meyer} \cite[Chap.\ IV, Th\'eor\`eme 30]{DM}),
which we continue to denote by $v$. We note that thanks to \eqref{rd:momentb}, $v\in L^p_{\text{\it loc}}(dtdx)$ a.s.

   We apply Proposition \ref{pr:Garsia} to this version of $v$ with
\begin{equation}
   S=I\times J,\qquad \rho((t\,,x)\,;(s\,,y))= {\bf \Delta}
   ((t\,,x)\,;(s\,,y)), \qquad \mu(dt\, dx) = dtdx,
\end{equation}
and
\begin{equation}
   \Psi(x)=\vert x \vert^{p},\qquad \Psi^{-1}(y)=y^{1/p},
   \qquad p(x) = \vert x \vert^{\alpha + (6/p)}.
\end{equation}
Let
\begin{equation}
   \mathcal{C} = \int_S dtdx \int_S dsdy\
   \frac{\vert v(t\,,x) - v(s\,,y)\vert^p}{
   [{\bf \Delta}((t\,,x)\,;(s\,,y))]^{6+\alpha p}}.
\end{equation}
By (\ref{rd:momentb}),
\begin{equation}\begin{split}
   \E[\mathcal{C}] &\le \int_S dtdx \int_S dsdy\
   	[{\bf \Delta}((t\,,x)\,;(s\,,y))]^{q-3-\alpha p} \\
   &\le 4 \vert I\vert\, \vert J \vert \int_0^{\vert I\vert}
   	d\tilde{u} \int_0^{\vert J \vert} dv\ (\tilde{u}^{1/2} +v)^{q-3-\alpha p}.
\end{split}\end{equation}
We can check readily that the preceding integral is finite
using only the fact that $\alpha \in [0\,,q/p[$. Therefore, 
\begin{equation}\label{EC}
   \E[\mathcal{C}] <  \infty.
\end{equation}

   Since $v\in L^p_{\text{\it loc}}(dtdx)$ a.s., and because $p>1$, a well-known theorem of Jessen, Marcinkiewicz,
and Zygmund implies that the following
holds with probability one:
\begin{equation}
	\lim_{\e,\delta\downarrow 0}\frac{1}{4\e\delta}
	\int_{t-\e}^{t+\e} \int_{x-\delta}^{x+\delta} v(a\,,b)\, da\, db
	= v(t\,,x),
\end{equation}
for almost all $(t\,,x)\in I\times J$. See
\textsc{Khoshenvisan}
\cite[Theorem 2.2.1, Chapter 2, p.\ 58]{Khoshnevisan:02}. In particular,
\eqref{cond:leb} holds in the present setting.

   We now take into account Remark \ref{remgarsia}(b),
   and deduce that for a.a.\
   $\omega$ there exists a set $D(\omega) \subset S$
   with full Lebesgue measure such that for all
   $(t\,,x), (s\,,y) \in D(\omega)$,
	\begin{equation}\begin{split}
	   &\vert v(t\,,x)(\omega) - v(s\,,y)(\omega)\vert \\
	   &\hskip1in \le 8 \sup_{(r,\bar y)} \int_0^{2 {\bf \Delta}
	   	((t\,,x)\,;(s\,,y))} \Psi^{-1}\left(\frac{\mathcal{C}}{
			[\mu(B_{\rho}((r\,,\bar y)\,,u/2)]^2} \right)
			u^{\alpha - 1 + (6/p)}\, du.
	\end{split}\end{equation}
	One can check directly that there exists a $c>0$ such that
	$\mu(B_{\rho}((r\,,\bar y)\,,u/2) \geq cu^3$
	for all $u>0$ and
	$(r\,,\bar y) \in S$. Therefore,
	\begin{equation}\label{A20}\begin{split}
	   \vert v(t\,,x)(\omega) - v(s\,,y)(\omega)\vert &
	   	\le 8 \int_0^{2 {\bf \Delta}((t,x)\,;(s,y))}
			\mathcal{C}^{1/p} u^{\alpha -1} du\\
	   &= 2^\alpha\, 8\, \mathcal{C}^{1/p}\,
	   	[{\bf \Delta}((t\,,x)\,;(s\,,y))]^\alpha.
	\end{split}\end{equation}
	Define
	\begin{equation}
		\tilde{v}(t\,,x)(\omega):= 
		\limsup_{(s,y) \in D(\omega): \, (s,y)
		\rightarrow (t,x)} v(s\,,y)	(\omega).
	\end{equation}
	Since $v(\cdot)(\omega)$ is uniformly continuous on $D$,
	$\tilde{v}(\cdot)(\omega)$ is continuous on $\bar D(\omega) = S$
	and coincides with $v(\cdot)(\omega)$ in $D(\omega)$.
	In addition, by (\ref{rd:momentb}), $v(s\,,y)$ converges
	to $v(t\,,x)$ in $L^p$ as $(s\,,y)$ converges to $(t\,,x)$.
	Therefore, $\tilde{v}(t\,,x) = v(t\,,x)$ a.s.\ 
	for all $(t\,,x) \in S$, and hence
	$\tilde v$ is a continuous version of $v$. By (\ref{A20}), 
	\begin{equation}
	   \left(\sup_{(t,x) \neq (s,y)}\frac{\vert \tilde v(t\,,x)
	   - \tilde v(s\,,y)\vert}{[{\bf \Delta}
	   ((t\,,x)\,;(s\,,y))]^{\alpha}} \right)^p
	   \le 2^{\alpha p} 8^p \mathcal{C}.
	\end{equation}
	Equation (\ref{A.12}) now follows from (\ref{EC}).
\end{proof}

\section{Appendix: On Energy Reduction for Smoothed Measures}\label{apendA}

The goal of this appendix is to prove precise
versions of the statement,
``if we smooth a measure then we lower its
energy.''

\begin{theorem}\label{th:energybound:1}
	Let $0<\alpha<d$ and $\mu$ be a probability measure on $\R^d$.
	Then for all probability density functions $g:\R^d\to\R_+$ with compact support,
	\begin{equation}
		I_\alpha(g*\mu) \le I_\alpha(\mu).
	\end{equation}
\end{theorem}

\begin{theorem}\label{th:energybound:2}
	Choose and fix $n>1$.
	Then there exists a positive and finite
	constant $c$---depending only on
	$(d\,,n)$---such that for all probability measures $\mu$  on $[-n\,,n]^d$
	and all probability density functions $g:\R^d\to\R_+$ with compact support,
	\begin{equation}
		I_0(g*\mu) \le c\, I_0(\mu).
	\end{equation}
\end{theorem}
The proof requires some terminology
from harmonic analysis.
A function $\kappa:\R^d \to\R\cup\{\infty\}$ is called
a \emph{potential kernel} if: 
(i) $\kappa(x)\ge 0$ for all $x\neq 0$; 
(ii) $\kappa $ is continuous on $\R^d\setminus\{0\}$; and
(iii) $\kappa(0)=\infty$; $\kappa$ is called \emph{of positive
type} if its Fourier transform $\hat{\kappa}$ --- 
viewed in the sense of
distributions --- is a nonnegative function. We choose
the following normalization of Fourier transforms:
$\hat{\kappa}(\xi) = \int_{\R^d} \exp(i\xi\cdot x) \kappa(x)\, dx$
for all $\xi\in\R^d$ and $\kappa\in L^1(\R^d)$.

The following is well known; see for example 
\textsc{Kahane} \cite[Remark 2, p.\ 133]{Kahane:68}.

\begin{proposition}\label{pr:energy:fourier}
	If $\kappa:\R^d\to\R_+\cup\{\infty\}$ is a potential
	kernel of positive type, then for all Borel probability
	measures $\mu$ on $\R^d$,
	\begin{equation}
		\iint \kappa(x-y)\, \mu(dx)\, \mu(dy)
		=\frac{1}{(2\pi)^d} \int_{\R^d} \left| \hat\mu(\xi)
		\right|^2 \hat\kappa(\xi)\, d\xi.
	\end{equation}
\end{proposition}

First we prove Theorem \ref{th:energybound:1}; 
it is technically simpler than Theorem \ref{th:energybound:2},
and yet affords us the chance to discuss the reasons
for the veracity of both theorems.

\begin{proof}[Proof of Theorem \ref{th:energybound:1}]
	Define $\kappa(x) := \|x\|^{-\alpha}$, where
	$1/0:=\infty$, to find that
	$\kappa$ is a potential kernel of positive type
	with $\hat{\kappa}(\xi)=c\|\xi\|^{-d+\alpha}$;
	see \textsc{Stein} \cite[Chap.V, \S 1, Lemma 2(b)]{Stein}, or \textsc{Kahane} \cite[p.\ 134]{Kahane:68}, for
	example. Define $\nu(dx):=(g*\mu)(dx)$
	and apply Proposition \ref{pr:energy:fourier}
	with $\nu$ in place of $\mu$ to find that
	\begin{equation}
		I_\alpha(g*\mu) = \frac{1}{(2\pi)^d}
		\int_{\R^d} \left| \hat{g}(\xi)\right|^2
		\left|\hat\mu(\xi)\right|^2
		\hat{\kappa}(\xi)\, d\xi.
	\end{equation}
	Because $|\hat{g}(\xi)|\le 1$, another appeal to
	Proposition \ref{pr:energy:fourier} completes
	the proof.
\end{proof}

Now we work to prove the more difficult
Theorem \ref{th:energybound:2}.
We define a function $\rho:\R^d\to\R_+\cup\{\infty\}$ by
\begin{equation}
	\rho (x) := \frac{\exp(-\|x\|)}{\|x\|^{d/2}},
\end{equation}
where $\rho(0):=\infty$. 
Define $\kappa:\R^d\to\R_+\cup\{\infty\}$
by
\begin{equation}
	\kappa(x) := (\rho* \rho)(x)
	=\int_{\R^d} \rho(x-y)\, \rho(y)\, dy.
\end{equation}

\begin{lemma}\label{lem:kappa}
	The function $\kappa$ is an integrable potential kernel
	of positive type.
\end{lemma}

\begin{proof}
	Because $\rho(x)\ge 0$ is measurable
	the convolution is a well-defined 
	nonnegative Borel-measurable function
	on $\R^d$. Standard arguments show
	that $\kappa$ is at least as smooth as $\rho$.
	Since $\rho$ is continuous on $\R^d\setminus\{0\}$,
	then so is $\kappa$. Because $\kappa(0)$ is
	manifestly infinite, this proves that $\kappa$ is
	a potential kernel. We may note that 
	$\|\kappa\|_1=\|\rho\|_1^2<\infty$. Therefore,
	$\hat{\kappa}$ is the $L^1$-form of the Fourier
	transform of $\kappa$. Finally, since $\rho$ is
	even, $\hat\rho$ is real-valued, and therefore $\hat{\kappa}(\xi)=|\hat{\rho}(\xi)|^2\ge 0$.
The lemma follows.
\end{proof}

\begin{lemma}\label{lem:kappa2}
	Let $N_0$ be as in \textnormal{(\ref{k})}.
	Then there exist positive and finite constants
	$c_1$ and $c_2$---depending only on $(d\,,N_0)$---such
	that for all $x\in B(0\,,N_0/2)$,  
	\begin{equation}
		c_1 {\rm K}_0(\|x\|)\le \kappa(x) \le
		c_2 {\rm K}_0(\|x\|).
	\end{equation}
\end{lemma}

\begin{proof}
	Choose and fix $x$ with $0<\|x\|\le N_0/2$, and write
	\begin{equation}
		\kappa(x) = T_1 + T_2 + T_3,
	\end{equation}
	where
	\begin{equation}\begin{split}
		T_1 & := \int_{\|y\|< 2\|x\|} \rho(x-y)\rho(y)\, dy,\\
		T_2 & := \int_{2\|x\| \le \|y\|\le 10 N_0} \rho(x-y)\rho(y)\, dy,\\
		T_3 &:= \int_{\|y\|>10 N_0} \rho(x-y)\rho(y)\, dy.
	\end{split}\end{equation}
	We estimate each $T_i$ separately.
	
	It will turn out that the main contribution to $\kappa(x)$
	comes from $T_2$. Therefore, we begin by
	bounding that quantity: If $2\|x\| \le \|y\|$, then
	$\|x-y\| \le \frac{3}{2}\|y\|$; thus,
	\begin{equation}\begin{split}
		T_2 &\ge \left( \frac{2}{3}\right)^{d/2}
			\int_{2 \|x\|\le \|y\|\le 10 N_0}
			\frac{e^{-3\|y\|/2}}{\|y\|^{d/2}}\rho(y)\, dy\\
		&\ge C_1 \int_{2 \|x\| \le \|y\|\le 10 N_0}
			\frac{dy}{\|y\|^d}.
	\end{split}\end{equation}
	We integrate this in polar coordinates to find that
	$T_2 \ge C_2(\ln N_0 + \ln(1/\|x\|))$. Because $T_1,T_3\ge 0$,
	it follows that $\kappa(x)$ is bounded below by a
	constant multiple of $\ln(N_0/\|x\|)$. This proves
	half of the lemma.
	
	For the other half, we note that if $2\|x\| \le \|y\|$, then
	$\|x-y\|\ge \|y\|/2$. Therefore,
	we can use an argument, similar to the one we used to bound $T_2$
	from below, in order to prove that
	\begin{equation}
		T_2 \le C_3 (\ln(10 N_0) +\ln(1/\|x\|)),
	\end{equation}
and since $\|x\| \le N_0/2$, the right-hand side is bounded above by $C_4 (\ln N_0 +\ln(1/\|x\|))$, provided $C_4$ is chosen large enough.
	
Next we bound $T_3$. Note that if $\|y\|>10 N_0$, then
	$\|x-y\|\ge 9 N_0$. Consequently,
	$\rho(x-y)$ is bounded from above, and hence
	$T_3 \le C_4 \int_{\R^d} \rho(y)\, dy <\infty$.
	
	Finally, we estimate $T_1$ by first writing it as
	\begin{equation}
		T_1 = T_{11} + T_{12},
	\end{equation}
	where
	\begin{equation}\begin{split}
		T_{11} &:= \int_{\substack{\|y\|\le 2\|x\| \\
			\|y-x\|\ge \|x\|/2}} \rho(x-y)\rho(y)\, dy,\\
		T_{12} &:= \int_{\substack{\|y\|\le 2\|x\| \\
			\|y-x\|< \|x\|/2}} \rho(x-y)\rho(y)\, dy.
	\end{split}\end{equation}
	
	If $\|y-x\|\ge\|x\|/2$, then $\rho(x-y)\le 2^{d/2} \|x\|^{-d/2}$,
	and thus,
	\begin{equation}\begin{split}
		T_{11} &\le \frac{2^{d/2}}{\|x\|^{d/2}}
			\int_{\substack{\|y\|\le 2\|x\| \\
			\|y-x\|\ge \|x\|/2}} \frac{\exp(-\|y\|)}{\|y\|^{d/2}}\, dy\\
		&\le \frac{2^{d/2}}{\|x\|^{d/2}}
			\int_{\|y\|\le 2\|x\|}\frac{dy}{\|y\|^{d/2}}\\
		&\le C_5.
	\end{split}\end{equation}
	The last line follows from integrating in polar coordinates.
	
	In order to estimate the remaining term
	$T_{12}$, we note that if $\|y\|\le 2\|x\|$ and
	$\|y-x\|<\|x\|/2$, then $\|y\|\ge \|x\|/2$, and 
	hence $\rho(y)\le 2^{d/2}\|x\|^{-d/2}$. Consequently,
	\begin{equation}\begin{split}
		T_{12} &\le \frac{2^{d/2}}{\|x\|^{d/2}}
			\int_{\|y-x\|\le\|x\|/2} \rho(x-y)\, dy\\
		&\le \frac{2^{d/2}}{\|x\|^{d/2}}
		\int_{\|z\| \le \|x\|/2} \frac{dz}{\|z\|^{d/2}}\\
		&\le C_5,
	\end{split}\end{equation}
	for the same constant $C_5$ as before.
	These remarks together prove the lemma.
\end{proof}

Now we prove Theorem \ref{th:energybound:2}.

\begin{proof}[Proof of Theorem \ref{th:energybound:2}]
	Thanks to Lemma \ref{lem:kappa2},
	\begin{equation}
		I_0(g*\mu) \le \frac{1}{c_1}\iint \kappa(x-y)\, 
		\nu(dy)\, \nu(dx),
	\end{equation}
	where $\nu(dx) := (g*\mu)(x)\, dx$.
	Because $|\hat\nu(\xi)|=|\hat{g}(\xi)\hat{\mu}(\xi)|\le
	|\hat{\mu}(\xi)|$, 
	Lemma \ref{lem:kappa} and
	Proposition \ref{pr:energy:fourier} together imply that
	\begin{equation}\begin{split}
		I_0(g*\mu) &\le\frac{1}{c_1(2\pi)^d}
			\iint\left| \hat\mu(\xi)\right|^2 \hat\kappa(\xi)\, d\xi\\
		&= \frac{1}{c_1}\iint \kappa(x-y)\, \mu(dx)\, \mu(dy).
	\end{split}\end{equation}
	Another application of Lemma \ref{lem:kappa2} shows that
	the latter term is at most $(c_2/c_1) I_0(\mu)$,
	whence follows the theorem with $C:= c_2/c_1$.
\end{proof}

\noindent\textbf{Acknowledgements.} Peter Imkeller and Michael Scheutzow shared with us their notes on Garsia's lemma and its variants. We thank them both.

\end{document}